\documentstyle[graphicx]{article}[15pt]
\input amssym.def
\catcode`@=12
 \long\def\@makefntext#1{\noindent #1}
\newskip\tabcentering \tabcentering=1000pt plus 1000pt minus 1000pt

\def\MCH#1#2{\setbox0=\hbox{\raise#1\hbox{#2}}\smash{\box0}}

\def\@evenfoot{}\def\@oddfoot{}

\def\sec#1{\vspace{5mm}\leftline{\bf #1}\vspace{3mm}}

\catcode`@=12

\def\bc{\begin{center}}
\def\ec{\end{center}}

\def\hang{\hangindent\parindent}
\def\textindent#1{\indent\llap{\qquad #1\ \ \enspace}\ignorespaces}
\def\ref{\par\hang\textindent}

\def\a1{(a_1, a_2, \cdots, a_n)}

\def\a{\alpha}

\begin{document}
\thispagestyle{empty}
\vspace*{-3.0truecm}
\noindent
\vspace{1 true cm}
 \bc{\large\bf Concentration  in vanishing adiabatic exponent limit of  solutions to the  Aw-Rascle traffic model$^{ **}$

\footnotetext{$^{*}$Corresponding author. Tel: +86-0591-83852790.\\
\indent \,\,\,\,\,\,\,\,E-mail address:  zqshao@fzu.edu.cn.\\
\indent \,\,\,\,\,$^{**}$Supported by  the Scientific Research
Foundation of the
 Ministry of Education of China (No. 02JA790014),
  the Natural Science Foundation of Fujian Province of China   (No.
 2015J01014)   and  the Science and Technology Developmental Foundation
of Fuzhou University (No. 2004-XQ-16).  }}\ec

  \vspace*{0.2 true cm}
\bc{\bf Shouqiong Sheng$^{a}$, Zhiqiang  Shao$^{a, *}$   \\
{\it $^{a}$College of Mathematics and Computer Science, Fuzhou University, Fuzhou 350108, China}
 }\ec

 \vspace*{2.5 true mm}
\setlength{\unitlength}{1cm}
\begin{picture}(20,0.1)
\put(-0.6,0){\line(1,0){14.5}}
\end{picture}

 \vspace*{2.5 true mm}
\noindent{\small {\small\bf Abstract}

 \vspace*{2.5 true mm} In this paper, we study  the phenomenon of concentration and the
formation of delta shock wave in vanishing adiabatic exponent limit of  Riemann solutions to the  Aw-Rascle traffic model. It is proved that as the adiabatic exponent vanishes, the limit of solutions tends to a special delta-shock rather than the
classical one to the zero pressure gas dynamics. In order to further study this problem, we
consider a perturbed Aw-Rascle model and proceed to investigate the limits
of solutions. We rigorously proved that, as $\gamma$ tends to one, any
Riemann solution containing two shock waves tends to a delta-shock to  the zero pressure gas dynamics
 in the distribution sense. Moreover, some representative
numerical simulations are exhibited to confirm the theoretical analysis.

 \vspace*{2.5 true mm}
\noindent{\small {\small\bf MSC: } 35L65;  35L67

 \vspace*{2.5 true mm}
\noindent{\small {\small\bf Keywords:}  Aw-Rascle traffic model; Riemann solutions; Delta shock wave; Vanishing adiabatic exponent limit;
Zero pressure gas dynamics;
Weighted Dirac-measure;
Numerical simulations

 \vspace*{2.5 true mm}
\setlength{\unitlength}{1cm}
\begin{picture}(20,0.1)
\put(-0.6,0){\line(1,0){14.5}}
\end{picture}



\baselineskip 15pt
 \sec{\Large\bf 1.\quad  Introduction }
  The celebrated Aw-Rascle (AR)  model of traffic flow reads  (cf. [1]):
$$ \left\{\begin{array}{ll}\rho_{t}+(\rho u)_{x}=0,\\(\rho (u+p(\rho)))_{t}+
 (\rho u(u+p(\rho)) )_{x}=0,\end{array}\right .\eqno{(1.1)}
$$
where $\rho$      and $u$    represent the traffic density and velocity of the cars  located at position $x$ at time $t$, respectively; $p $  is the velocity offset and
called as the ``pressure" inspired from gas dynamics.
  The model (1.1) is now widely used to study the formation and dynamics of traffic jams. It was proposed by Aw and Rascle [1] to remedy the deficiencies of second order models of car traffic pointed out by Daganzo [6] and had also been independently derived by Zhang [30]. Since its introduction, it had received extensive attention (see [18, 20, 23, 28], etc.).

 In this paper, we are concerned with  the ``pressure"
function  $$p(\rho)=\rho^{\gamma}, \,\,\,\,\,0<\gamma< 1.\eqno{(1.2)}$$

 The  Riemann solutions of (1.1) with classical pressure $p(\rho)=\rho^{\gamma}$  ($\gamma>0$)
were  obtained at low densities by Aw and Rascle [1]. Lebacque, Mammar, and Salem1 [13] also solved
the Riemann problem of (1.1) with classical pressure $p(\rho)=\rho^{\gamma}$  ($\gamma>0$) with an extended fundamental diagram
 for all possible initial data. Sun [28] studied the interactions of elementary waves
to  system (1.1).

We are interested in
   the Riemann problem for (1.1)-(1.2)   with initial data
$$ (\rho, u)(0, x) =\left\{\begin{array}{ll} (\rho_{-},
u_{-}),\,\,\,\,x< 0,\\(\rho_{+},
u_{+}),\,\,\,\,x> 0,\end{array} \right.\eqno{(1.3)}$$
where  $\rho_{\pm}>0$ and $u_{\pm}$   are given constant states. We assume that
$u_{+}<u_{-}$.

  System (1.1)-(1.2) is just like a hyperbolic system for conservation laws of the form
  $$ \partial_{t}U+\partial_{x}F(U)=0,\eqno{(1.4)} $$
with
$$ U=\left(\begin{array}{c}\rho \\\rho u+\rho^{\gamma+1}
 \end{array}\right), ~~F(U)=\left(\begin{array}{c}\rho u \\\rho u^{2}+u\rho^{\gamma+1}
 \end{array}\right)=0.$$

When $\gamma\rightarrow 0$, the limiting system of (1.1)-(1.2) formally becomes   the zero pressure gas dynamics,
$$ \left\{\begin{array}{ll} \rho_{t}+(\rho u)_{x}=0,\\(\rho u)_{t}+
 (\rho u^{2} )_{x}=0,\end{array}\right .\eqno{(1.5)}
$$
 which can be used to describe the process of the motion of free particles sticking under collision
and depict the formation of large scale in the universe.  The solutions to the zero pressure gas dynamics were widely studied by many scholars (see [2-3, 7-9, 15-16, 26], etc. ). In particular, the   existence of measure solutions
of the Riemann problem was first proved by Bouchut [2] and the existence of the global weak solution
was obtained by Brenier and Grenier [3] and E, Rykov and Sinai [7]. Sheng and Zhang [26] discovered
that the $\delta$-shocks and vacuum states do occur in the Riemann solutions to the zero pressure gas dynamics (1.5) by
the vanishing viscosity method. Huang and Wang [9] proved the uniqueness of the weak solution
for the case when the initial data is a Radon measure.

 A distinctive feature for (1.5) is just that the $\delta$-shocks and vacuum states do occur in the Riemann solutions.  In
paper [23], Shen and Sun studied the limits of Riemann solutions of (1.1)  with classical pressure $ p(\rho) =\varepsilon\rho^{\gamma}$
$(\gamma> 0)$  as $\varepsilon \rightarrow 0+$. They
identified a special $\delta$-shock  in the limit of solutions, whose the propagation
speed and the strength
are different from  those of the zero pressure gas dynamics  (1.5).   Then,
they analyzed a perturbed Aw-Rascle model and proved that the limit of Riemann solutions to the
perturbed Aw-Rascle model are those of (1.5) when $\varepsilon \rightarrow 0+$. The idea of vanishing pressure limits dates back  to early works of   Li [14], Chen and Liu [4,5], and the vanishing pressure limit
method was also applied to other systems [17-20, 22, 24-25, 29].

Let us turn to the Euler system  of power law in Eulerian coordinates,
$$ \left\{\begin{array}{ll} \rho_{t}+(\rho u)_{x}=0,\\(\rho u)_{t}+
 (\rho u^{2}+p(\rho) )_{x}=0,\end{array}\right .
\eqno{(1.6)}$$  When the pressure tends to zero or a constant, the Euler system  (1.6)  formally tends to the zero pressure gas dynamics.  In earlier seminal papers,
Chen and Liu   [4] first showed the formation of $\delta$-shocks and vacuum states of the Riemann solutions to the Euler system (1.6) for polytropic gas by taking limit $\varepsilon \rightarrow 0+$ in the model $p(\rho)=\varepsilon\rho^{\gamma}/\gamma$   $ (\gamma >1)$,
which describe the phenomenon of concentration and cavitation rigorously in mathematics.  Further, they also obtained the same results for the Euler
equations for nonisentropic
fuids in [5]. The same problem for the  Euler equations  (1.6) for
isothermal case $( \gamma =1)$ was studied by Li [14].  Recently,  Muhammad Ibrahim, Fujun Liu and Song Liu
[10] showed the same phenomenon of concentration also exists in the mode $p(\rho) = \rho^{\gamma}$ $(0 < \gamma < 1)$ as $\gamma \rightarrow 0$, which is the case that the pressure goes to a constant. Namely, they  showed rigorously the formation of delta wave with the limiting behavior of Riemann solutions to the Euler equations (1.6).

Motivated by [10], for the Aw-Rascle model (1.1) with classical pressure (1.2), we show the same phenomenon of concentration also exists in the case $0 < \gamma < 1$ and  $u_{+}<u_{-}$ as   $\gamma \rightarrow 0$. We can see that,  as $\gamma\rightarrow 0$, the Riemann solution converges to a special delta shock solution, whose the propagation speed and the strength are different from those of the PGD model (1.5), which means the Riemann solution of (1.1)-(1.2) don't converge to the delta shock solution of (1.5).

In order to solve this problem, we motivated by [23], adding a suitable perturbation in the pressure term in the Aw-Rascle model (1.1)-(1.2). That is we consider the perturbed Aw-Rascle (PAR) model as follows:
$$ \left\{\begin{array}{ll}\rho_{t}+(\rho u)_{x}=0,\\
\bigg(\rho u+\frac{1}{\gamma}\rho^{\gamma}\bigg)_{t}+\bigg(\rho u^{2}+u\rho^{\gamma} \bigg)_{x}=0,\end{array}\right .\eqno{(1.7)}
$$
where $1<\gamma<3$.  For convenience and conciseness, we replace $\rho p(\rho)$ with $p(\rho)$ in (1.1) and take
$p(\rho) = \rho^{\gamma}$ for $\gamma\in(1,3)$. In the system (1.7), $p(\rho) = \rho^{\gamma}$ can be regarded as the traffic pressure term and  $1<\gamma<3$ is analogous with the adiabatic exponent
 $0<\gamma<2$ in the Aw-Rascle model (1.1)-(1.2).  It
is proved that when $\gamma\rightarrow 1$, the limit of the Riemann solutions containing two shock waves of the perturbed
Aw-Rascle model is exactly a delta shock solution of
 the zero pressure gas dynamics (1.5).

Finally, by using the fifth-order
weighted essentially non-oscillatory scheme and third-order Runge-Kutta method  [12, 27], some representative numerical simulations are
exhibited, which are completely consistent with theoretical analysis.

\vskip 0.1in
The rest of the paper is organized as follows. For the sake of completeness, in Section 2, we briefly
 review the delta shock wave and vacuum state in the Riemann solutions
of
 the zero pressure gas dynamics
(1.5). In Section 3, we display some results on  the Riemann solutions of (1.1)-(1.2)  when $0 < \gamma < 1$. In Section 4, we discuss  the limits of
Riemann solutions of  (1.1)-(1.2) as the adiabatic exponent vanishes. In Section 5, we display some results on  the Riemann solutions of (1.7)  when $1 < \gamma < 3$. In Section 6,   we show rigorously the formation of delta shock wave
with the limiting behavior of Riemann solutions  of (1.7) as  $\gamma\rightarrow 1$.
In section 7, we present the numerical results.

\baselineskip 15pt
 \sec{\Large\bf 2.\quad   Preliminaries }
  For the sake of completeness, in this section we briefly recall  the delta shock wave and vacuum state  in the Riemann solutions of
 the zero pressure
gas dynamics (1.5). More details can be found in [26, 24, 16, 11].

The system (1.5)
 has  a double
eigenvalue
$ \lambda=u
$
and only one  right eigenvector
$ \overrightarrow{r}=(1, 0)^{T}.
$
The system is obviously nonstrictly hyperbolic,  and  $\lambda$  is
 linearly degenerate  by
$ \nabla\lambda\cdot \overrightarrow{r}\equiv 0,
   $,  in which $\bigtriangledown$  denotes the gradient with respect to $(\rho,u)$. Therefore, in classical sense, the associated  elementary waves involve
only contact discontinuities.  It can
be seen from previous works [11,16, 24, 26] that
the Riemann problem for (1.5) with initial data (1.3) can be solved by contact discontinuities, vacuum  or delta shock wave connecting two constant states $(\rho_{\pm},u_{\pm})$.

When $u_{-}< u_{+}$, there is no characteristic passing through the region $u_{-}t < x < u_{+}t $  and  the vacuum  appears in this  region. The solution can be expressed as
$$ (\rho,u)(t,x)=\left\{\begin{array}{ll} (\rho_{-}, u_{-}),\,\,\,\,\,\,-\infty<x<u_{-}t, \\(0,\frac{x}{t}),\,\,\,\,\,\,\,\,\,\,\,\,\,\,\,\,\,u_{-}t \leq x \leq u_{+}t,\\(\rho_{+},u_{+}), \,\,\,\,\,\,u_{+}t<x < +\infty.
\end{array}\right .\eqno{(2.1)}
$$

When $u_{-}= u_{+}$,  the constant states $(\rho_{\pm},u_{\pm})$ can be connected by a contact discontinuity. The solution can be expressed as
$$(\rho, u)(t,x)=\left\{
    \begin{array}{ll}
      (\rho_{-}, u_{-}), \,\, \,\,\,\,\,-\infty<x<u_{-}t,\\
     (\rho_{+}, u_{+}),\,\,\,\,\,\,\,u_{-}t<x<+\infty.
    \end{array}
  \right.\eqno{(2.2)}$$

When $u_{-}> u_{+}$,   the characteristic lines from initial data  will overlap,  so the Riemann solution cannot be constructed by using the classical waves,
we seek   a  solution containing a weighted Dirac delta function with the support on a line.

To do so, a
  two-dimensional weighted delta function  $w(s)\delta_{S}$ supported on a smooth curve $ S=\{(t(s),x(s)):a<s<b\}$ is defined by
$$ \langle w(t)\delta_{S},\varphi(t,x)\rangle=\int_{a}^{b}w(t(s))\varphi(t(s),x(s))ds, \eqno{(2.3)}$$
for all test functions $\varphi(t,x) \in C_{0}^{\infty}([0,+\infty)\times (-\infty, +\infty)).$ \\

For the Riemann problem  with $u_{+}<u_{-}$, we can construct a dirac-measured solution with
 parameter $\sigma$  as follows,
 $$  \rho(t,x)=\rho_{0}(t,x)+w(t)\delta_{S}, ~~u(t,x)=u_{0}(t,x), \eqno{(2.4)}
 $$
 where $ S=\{(t,\sigma t):0\leq t<+\infty\} $,
 $$\rho_{0}(t,x)=\left\{
                  \begin{array}{ll}
                    \rho_{-}, & \hbox{$x<\sigma t $,} \\
                    \rho_{+}, & \hbox{$x>\sigma t$,}
                  \end{array}
                \right.
\eqno{(2.5)}$$$$u_{0}(t,x)=\left\{
                  \begin{array}{ll}
                   u_{-}, & \hbox{$x<\sigma t $,} \\
                    \sigma, & \hbox{$x=\sigma t$,} \\
                    u_{+}, & \hbox{$x>\sigma t$,}
                  \end{array}
                \right.
\eqno{(2.6)}$$
and$$ w(t)=t(\sigma[\rho]-[\rho u]),\eqno{(2.7)}$$
in which $[q]=q_{+} -q_{-}$ denotes the jump of  function $q$  across the discontinuity discontinuity.
The dirac-measured
solution $(\rho, u)$ constructed above is called a delta shock solution of (1.5)
in the sense of distributions  if
$$ \langle \rho,\varphi_{t}\rangle+\langle \rho u,\varphi_{x}\rangle=0,\eqno {(2.8)}
$$$$ \langle \rho u,\varphi_{t}\rangle+\langle \rho u^{2}, \varphi_{x}\rangle=0,\eqno {(2.9)}
$$
hold  for any test function $\varphi (t,x)\in C_{0}^{\infty}([0,+\infty)\times (-\infty, +\infty))$, where
$$ \langle \rho, \varphi\rangle=
       \int_{0}^{+\infty}\int_{-\infty}^{+\infty}\rho_{0}(t,x)\varphi (t,x)dxdt+\langle w(t)\delta_{S}, \varphi(t,x)\rangle,$$
       $$\langle \rho u, \varphi\rangle=
      \int_{0}^{+\infty}\int_{-\infty}^{+\infty}\rho_{0}(t,x)u_{0}(t,x)\varphi (t,x)dxdt+\langle \sigma  w(t)\delta_{S}, \varphi(t,x)\rangle.
  $$
  Then the following generalized Rankine-Hugoniot relation
$$
\left\{
     \begin{array}{ll}
       \frac{dx}{dt}=\sigma, \\
       \frac{dw(t)}{dt}=\sigma [\rho]-[\rho u], \\
       \frac{d(w(t)\sigma)}{dt}=\sigma [\rho u]-[\rho u^{2}],
     \end{array}
   \right.\eqno (2.10)
$$
holds, where $[\rho]= \rho_{+}-\rho_{-}$,   with initial data
$$(x, w)(0) = (0,  0).\eqno{(2.11)}$$

\indent
To guarantee uniqueness,  the delta  shock  should satisfy the entropy
condition:
$$u_{+}<\sigma<u_{-},\eqno{(2.12)}$$
which means that all the characteristic lines on both sides of the discontinuity are incoming. So it is a overcompressive condition.

 Solving (2.10) with initial data
 (2.11) under the entropy condition (2.12),
we have$$
 w(t)=\sqrt{\rho_{-}\rho_{+}}(u_{-}-u_{+})\,t, ~~~~\sigma=\frac{\sqrt{\rho_{+}}u_{+}+\sqrt{\rho_{-}}u_{-}}{\sqrt{\rho_{+}}+\sqrt{\rho_{-}}}.\eqno{(2.13)}$$

Therefore, a delta shock solution defined by (2.4) with (2.5), (2.6) and (2.13) is obtained.

\baselineskip 15pt
 \sec{\Large\bf 3.\quad  Riemann solutions of the AR model (1.1)-(1.2)}In this section, we  review the Riemann solutions of  (1.1)-(1.2) with initial data (1.3), for which the detailed investigations can be found in Sun [28].

The system  (1.1)-(1.2) has two eigenvalues
 $$ \lambda_{1}=u-\gamma\rho^{\gamma},\,\,\,\,\,\,\,\lambda_{2}=u,\eqno{(3.1)}
$$
with the corresponding right eigenvectors
$$ \overrightarrow{r}_{1}=(1,-\gamma\rho^{\gamma-1})^{T},\,\,\,\,\,\,\overrightarrow{r}_{2}=(1, 0)^{T}\eqno{}
$$
satisfying$$\nabla\lambda_{1}\cdot \overrightarrow{r}_{1}=-\gamma(\gamma+1)\rho^{\gamma-1}<0,$$
and $$\nabla\lambda_{2}\cdot \overrightarrow{r}_{2}\equiv 0.\eqno{}
   $$
Therefore,  system  (1.1)-(1.2)  is strictly
hyperbolic for $\rho> 0$, and $\lambda_{1}$ is genuinely
nonlinear for $\rho> 0$ and the associated wave
is either shock wave or rarefaction wave, while $\lambda_{2}$ is always linearly degenerate and the associated wave is the contact  discontinuity.

Since  (1.1), (1.2) and the Riemann data (1.3) are invariant under stretching of coordinates:
$(t, x)\rightarrow (\tau t,\tau x)~(\tau$ is constant),  we seek the self-similar solution $$(\rho,u)(t,x)=(\rho,u)(\xi),\,\,\,\,\xi=\frac{x}{t}.$$
 Then the Riemann problem (1.1), (1.2) and (1.3) is reduced to the following boundary value problem of the ordinary differential equations:
$$ \left\{\begin{array}{ll}
   -\xi\rho_{\xi}+(\rho u)_{\xi}=0,\\
   -\xi(\rho u+\rho^{\gamma+1})_{\xi}+(\rho u^{2}+u\rho^{\gamma+1})_{\xi}=0,\end{array}\right .\eqno{(3.2)}
$$
with $(\rho,u)(\pm\infty)=(\rho_{\pm},u_{\pm}).$\\
\indent
For any smooth solution, system (3.2) can be written as
$$
\left(\begin{array}{cccc}u-\xi &\rho\\(u-\xi)( u+(\gamma+1)\rho^{\gamma}) &-\xi \rho+2\rho u+\rho^{\gamma+1}
 \end{array}\right)\left(\begin{array}{cccc}\rho_{\xi}\\u_{\xi}
 \end{array}\right)=0.\eqno{(3.3)}$$
Besides the constant solution $$(\rho,u)(\xi)={\mathrm constant }   \,\,\,\,\,(\rho> 0),$$ it provides a rarefaction wave which is
a continuous solution of (3.3) in the form $(\rho, u)(\xi)$. Then, for a given left state $(\rho_{-},u_{-})$, the
rarefaction wave curves in the phase plane, which are the sets of states that can be
connected on the right by a 1-rarefaction wave, are as follows:
$$R_{}(\rho_{-},u_{-}):
\left\{
  \begin{array}{ll}
    \xi=\lambda_{1}=u-\gamma\rho^{\gamma}, \\
    u-u_{-}=-(\rho^{\gamma}-\rho_{-}^{\gamma}), \\
    \rho<\rho_{-},u>u_{-}.
  \end{array}
\right.
\eqno{(3.4)}  $$
Differentiating
 the second equation of  (3.4) with respect to $\rho$ yields $$u_{\rho} = -\gamma\rho^{\gamma-1}<0,$$
 and
$$u_{\rho\rho} =-\gamma(\gamma-1)\rho^{\gamma-2}>0,$$
which mean that for $0<\gamma<1$, the rarefaction wave curve $R_{}(\rho_{-},u_{-})$ is monotonic decreasing and convex  in the $(\rho, u)$ phase plane
$(\rho> 0)$.  Moreover, it  can be concluded  from (3.4) that
 $\lim\limits_{\rho\rightarrow 0^{+}}u= u_{-}+\rho_{-}^{\gamma}$  for the rarefaction wave curve $R_{}(\rho_{-},u_{-})$, which implies that
$R_{}(\rho_{-},u_{-})$ intersects the $u$-axis at the point $(0,  \widetilde{u}_{\ast})$,  where $\widetilde{u}_{\ast}$ is determined by
$\widetilde{u}_{\ast}= u_{-}+\rho_{-}^{\gamma}$.

  \indent
For a bounded discontinuity at $\xi=\sigma,$ the Rankine-Hugoniot relation
$$\left\{
    \begin{array}{ll}
      -\sigma[\rho]+[\rho u]=0, \\
     -\sigma[\rho u+\rho^{\gamma+1}]+[\rho u^{2}+u\rho^{\gamma+1}]=0,
    \end{array}
  \right.\eqno (3.5)
$$
holds, where $[\rho]=\rho -\rho_{-},$   etc.  Eliminating $\sigma$ from (3.5), we obtain
$$[\rho][\rho u^{2}]-([\rho u])^{2}
     =-[\rho][u\rho^{\gamma+1}]+[\rho  u][\rho^{\gamma+1}].\eqno {(3.6)}
     $$
    Simplifying (3.6) yields $$ (u-u_{-})^{2}=-(u-u_{-})(\rho^{\gamma}-\rho_{-}^{\gamma}).
     \eqno $$If $u-u_{-}\neq 0$, we have
$$ u-u_{-}=-(\rho^{\gamma}-\rho_{-}^{\gamma})\, \,\, \mathrm{ and } \,\,\,\sigma=u-\frac{\rho_{-}(\rho^{\gamma}-\rho_{-}^{\gamma})}{\rho-\rho_{-}},
     \eqno {(3.7)}$$
where $\sigma$, $(\rho_{-},u_{-})$ and $(\rho_{},u_{})$ are the shock speed, the left state and the right state, respectively.

Otherwise, for case $u = u_{-}$ (i.e., $[u] = 0)$, we have
$$\sigma = u = u_{-}.$$

The classical Lax
 entropy conditions   imply  that the propagation speed $\sigma$ for the 1-shock wave
has to be satisfied with$$\sigma<\lambda_{1}(\rho_{-},u_{-}),\,\,\,\,\lambda_{1}(\rho,u)<\sigma<\lambda_{2}(\rho,u).
$$
From the first equation of (3.5), we obtain
$$\sigma=\frac{\rho u -\rho_{-}u_{-}}{\rho-\rho_{-}} = u_{-}+\frac{\rho}{\rho-\rho_{-}}(u-u_{-}).$$
If $u>u_{-}$, then from (3.7),  we have $\rho<\rho_{-}$, and
 $$\sigma- u_{-}=\frac{\rho}{\rho-\rho_{-}}(u-u_{-})=-\frac{\rho(\rho^{\gamma}-\rho_{-}^{\gamma})}{\rho-\rho_{-}}=-\rho\gamma
 \overline{\rho}^{\gamma-1},
$$
for some $\bar{\rho}\in(\rho,\rho_{-}).$  By direct calculation, we have
$$\gamma\rho_{-}^{\gamma}-\rho\gamma
 \overline{\rho}^{\gamma-1}>\gamma(\rho_{-}^{\gamma}-
 \rho^{\gamma})>0,
$$ which implies that
$$\sigma-u_{-} >-\gamma\rho_{-}^{\gamma}.
$$
This  contradicts with $\sigma<\lambda_{1}(\rho_{-},u_{-})$.
Then, given a left state $(\rho_{-},u_{-}),$
the possible states that can be connected to $(\rho_{-},u_{-})$ on the right by shock wave in the 1-family are as
follows:
$$S_{}(\rho_{-},u_{-}): \left\{
             \begin{array}{ll}
               \sigma=u-\frac{\rho_{-}(\rho^{\gamma}-\rho_{-}^{\gamma})}{\rho-\rho_{-}}, \\
             u-u_{-}=-(\rho^{\gamma}-\rho_{-}^{\gamma}), \\
              \rho>\rho_{-},u<u_{-}.
             \end{array}
           \right.
\eqno{(3.8)}  $$
Differentiating $u$ with respect to $\rho$ in the second equation of (3.8) gives
 that for $\rho>\rho_{-}$,
$$u_{\rho}=-\gamma\rho^{\gamma-1}<0\,\,\,\mathrm{and}\,\,\,u_{\rho\rho}=-\gamma(\gamma-1)\rho^{\gamma-2}>0,$$
which means that the shock wave curve $S_{}(\rho_{-},u_{-})$ is monotonic decreasing and convex in the $(\rho, u) $ phase plane  ($\rho>\rho_{-})$.
 It can also be derived from (3.8) that $\lim\limits_{\rho\rightarrow +\infty}u= -\infty$
 for the shock wave curve $S_{}(\rho_{-},u_{-})$, which indicates that the shock wave curve
intersects with the $\rho$-axis at a point.

Since $\lambda_{2}$ is linearly degenerate, the set of states $(\rho, u) $ can be connected to a given left state $(\rho_{-},u_{-})$
 by a contact discontinuity on the right if and only if
$$J: \xi=u=u_{-}. \eqno{(3.9)}$$

In the $(\rho, u) $ phase plane $( \rho, u \geq¡Ý 0)$, through a given point $(\rho_{-}, u_{-})$, we draw the elementary wave curves.  We find that the elementary wave curves divide the quarter phase plane $( \rho, u \geq¡Ý 0)$ into three regions,
$I=\{(\rho, u)|u<u_{-}\}$, $II=\{(\rho, u)|u_{-}<u<\widetilde{u}_{\ast}\}$,  and $III=\{(\rho, u)|u>\widetilde{u}_{\ast}\}$, where $\widetilde{u}_{\ast}= u_{-}+\rho_{-}^{\gamma}$, see Fig. 1.  According to the right
state $(\rho_{+},u_{+})$ in the different regions, one can construct  the
unique global Riemann solution connecting two constant states $(\rho_{\pm},u_{\pm})$ as follows:
(1) $(\rho_{+},u_{+})\in I(\rho_{-},u_{-}):$ $S+J,$
 (2)$(\rho_{+},u_{+})\in II(\rho_{-},u_{-}):$ $R+J,$
 (3)$(\rho_{+},u_{+})\in III(\rho_{-},u_{-}):$ $R+\mathrm{Vac}+J$
(see Fig. 1), where   ``+" means ``followed by".

\hspace{65mm}\setlength{\unitlength}{0.8mm}\begin{picture}(80,66)
\put(-24,2){\vector(0,2){50}}\put(-24,0){\line(0,4){52}}
 \put(4,0){\line(0,4){52}}\put(-48,0){\vector(2,0){103}} \put(3,-4){$u_{-}$} \put(-30,49){$\rho$}
\put(56,-1){$u$}
\put(6,45){$J$}\put(37,45){$J$}
\put(-18,35){$S_{}$}\put(-25,-4){$0$}\put(35,-4){$\widetilde{u}_{\ast}$}\put(35,0){\line(0,4){52}}
\put(5,19){$(\rho_{-},
u_{-})$}
\put(43,29){III }\put(24,6){$R_{}$ }\put(-5,29){I}
\put(17,29){II }\qbezier(35,0)(-14,12)(-24,42)
\end{picture}}
\vspace{0.6mm}  \vskip 0.2in \centerline{\bf Fig. 1.\, $(\rho, u)$-plane.
   } \vskip 0.1in \indent

\baselineskip 15pt
 \sec{\Large\bf 4.\quad   Limit of Riemann solutions  of the AR model (1.1)-(1.2) }In this section, we study the limiting behavior of the Riemann solutions of  (1.1)-(1.2) with the assumption $u_{+}<u_{-}$
as $\gamma$ tends to zero, that is, the formation of delta shock  as
 $\gamma\rightarrow0$  in the case $u_{+}<u_{-}$.

 \vskip 0.1in

\baselineskip 15pt
 \sec{\Large\bf 4.1.\quad   Formation of delta shock wave}

 For any fixed $\gamma \in(0, 1)$, when $u_{+}<u_{-}$,  namely $(\rho_{+}, u_{+})\in I( \rho_{-}, u_{-})$, the Riemann solution of
 (1.1)-(1.2)  is  a shock wave $S$ followed by a
contact discontinuity $J$ with the intermediate state $(\rho_{\ast},
u_{\ast})$
besides two constant states $( \rho_{-}, u_{-})$  and $(\rho_{+}, u_{+})$. They satisfy
$$S_{}: \left\{
             \begin{array}{ll}
               \sigma_{1}=u_{\ast}-\frac{\rho_{-}(\rho_{\ast}^{\gamma}-\rho_{-}^{\gamma})}{\rho_{\ast}-\rho_{-}}, \\
             u_{\ast}-u_{-}=-(\rho_{\ast}^{\gamma}-\rho_{-}^{\gamma}),\,\,\,
              \rho_{\ast}>\rho_{-},
             \end{array}
           \right.
\eqno{(4.1)}  $$
and
$$ J:\,\,\, \sigma_{2}=u_{\ast}=u_{+}, \,\,\,\,   \rho_{\ast}>\rho_{+},\eqno{(4.2)}
$$
where $\sigma_{1}$ and $\sigma_{2}$ are the propagation speeds of $S$ and $J$, respectively. Then we have the following lemmas.

 \vskip 0.1in
\noindent{\small {\small\bf Lemma 4.1.}
$\lim\limits_{\gamma\rightarrow0}\rho_{\ast}=+\infty,$ and $ \lim\limits_{\gamma\rightarrow 0}\rho_{\ast}^{\gamma}=:a
=1+u_{-}-u_{+}$.

\vskip 0.1in
\noindent{\small {\small\bf Proof.}   It follows from (4.1) and
(4.2) that
$$ u_{-}-u_{+}=\rho_{\ast}^{\gamma}-\rho_{-}^{\gamma}, \,\,\,\,   \rho_{\ast}>\rho_{\pm}.\eqno{(4.3)}
$$ Let $ \lim\limits_{\gamma\rightarrow0}\inf\rho_{\ast}=\alpha$, and $\lim\limits_{\gamma\rightarrow0}\sup\rho_{\ast}=\beta$.

If $ \alpha<\beta$ , then by the continuity of $\rho_{\ast}(\gamma)$, there exists a sequence $  \{\gamma_{k}\}_{k=1}^{\infty}\subseteq(0,1)$
such that
$$ \lim_{k\rightarrow +\infty}\gamma_{k}=0,\,\,\mathrm{ and}\,\, \lim_{k\rightarrow +\infty}\rho_{\ast}(\gamma_{k})=c,$$
for  some $ c\in(\alpha,\beta).$ Then substituting the sequence  into the right hand side of (4.3),  and  taking the limit $k\rightarrow +\infty$,
 we have
$$  u_{-}-u_{+}=\lim_{k\rightarrow+\infty}(\rho_{\ast}(\gamma_{k})^{\gamma_{k}}-\rho_{-}^{\gamma_{k}})=0.
\eqno{(4.4)} $$
This contradicts with the assumption  $u_{-}>u_{+}$.
Then we must have $\alpha=\beta$, which  means $\lim\limits_{\gamma\rightarrow1}\rho_{\ast}(\gamma)=\alpha.$

If  $\alpha\in(0,+\infty),$ then  we  can also get a contradiction when taking limit in (4.3). Hence $\alpha=0 $ or $ \alpha=+\infty$. By the condition
$\rho_{\ast}>\max\{\rho_{-},\rho_{+}\}$, it is easy to see that $\lim\limits_{\gamma\rightarrow 0}\rho_{\ast}(\gamma)=\alpha=+\infty.$

Next  taking the limit $ \gamma\rightarrow 0$ in (4.3), we have
$$u_{-}-u_{+}=\lim_{\gamma\rightarrow 0}(\rho_{\ast}^{\gamma}-\rho_{-}^{\gamma})=:a-1,$$
from which we  can get
$ a=1+u_{-}-u_{+}.$
 The proof is completed. $~\Box$

\vskip 0.1in
\noindent{\small {\small\bf Lemma 4.2.}

$$ \lim_{\gamma\rightarrow 0}\sigma_{1}=\lim_{\gamma\rightarrow 0}\sigma_{2}=\lim_{\gamma\rightarrow0}u_{\ast}=\sigma,\eqno{}
$$
where $\sigma=u_{+}$.

\vskip 0.1in
\noindent{\small {\small\bf Proof.} From
 (4.1), (4.2) and  Lemma 4.1, we immediately get
$$\lim_{\gamma\rightarrow0}\sigma_{1}=\lim_{\gamma\rightarrow0}\sigma_{2}=\lim_{\gamma\rightarrow0}u_{\ast}
=u_{-}-\lim_{\gamma\rightarrow0}(\rho_{\ast}^{\gamma}-\rho_{-}^{\gamma})
=u_{-}-(a-1)=u_{-}-
(u_{-}-u_{+})=u_{+}.
$$The proof is completed. $~~\Box$
\vskip 0.1in
Lemmas 4.1-4.2 show that when $\gamma$ tends to zero, $S$ and $J$ coincide, the intermediate
density $\rho_{\ast}$
 becomes singular.

\vskip 0.1in
\noindent{\small {\small\bf Lemma 4.3.}
$$ \lim\limits_{\gamma\rightarrow 0}\int_{\sigma_{1}}^{\sigma_{2}}\rho_{*}d\xi =\rho_{-}(u_{-}-u_{+})\neq 0. \eqno{(4.5)}
$$
\vskip 0.1in
\noindent{\small {\small\bf Proof.}
From the first equations of the Rankine-Hugoniot relation (3.5) for $S$ and $J$,   we have
 $$ \sigma_{1}(\rho_{-}-\rho_{\ast})=\rho_{-}u_{-}-\rho_{\ast}u_{\ast},\eqno{(4.6)}$$and $$\sigma_{2}(\rho_{\ast}-\rho_{+})=\rho_{\ast}u_{\ast}-\rho_{+}u_{+}. \eqno{(4.7)}$$
By  (4.6)
+(4.7), we get
 $$ \lim_{\gamma\rightarrow0}\rho_{\ast}(\sigma_{2}-\sigma_{1})=\lim_{\gamma\rightarrow0}(\rho_{-}u_{-}-\sigma_{1}\rho_{-}+\sigma_{2}\rho_{+}-\rho_{+}u_{+})
=\rho_{-}(u_{-}-u_{+}),
$$
which implies that$$\lim\limits_{\gamma\rightarrow 0}\int_{\sigma_{1}}^{\sigma_{2}}\rho_{*}d\xi =\rho_{-}(u_{-}-u_{+}).\eqno{(4.8)}$$
 The proof is completed. $~~\Box$
 \vskip 0.1in
 Lemma 4.3 shows that when $\gamma\rightarrow0$, the limit of $\rho_{\ast}$ has the same singularity as a weighted Dirac delta function at $\xi=u_{+}$.

\vskip 0.1in
\noindent{\small {\small\bf Remark 4.1.} It can be concluded from Lemmas 4.1-4.3 that, when $\gamma\rightarrow0$,   $S$ and $J$ coincide to form a new type of nonlinear hyperbolic wave, which is called as the delta shock wave in [45]. Compared with the  Riemann solutions
of (1.5), it is clear to see that the propagation speed and strength of the delta shock wave here are $\sigma = u_{+}$  and $w(t)=\rho_{-}(u_{-}-u_{+})\,t,$ which are
different from those of the
classical one to  the zero pressure gas dynamics (1.5).

\vskip 0.1in

Now, we give the following theorem which give a very nice depiction of the limit of Riemann solutions of (1.1) and (1.2) as $\gamma\rightarrow 0$  in the case $u_{+}<u_{-}.$

\vskip 0.1in
\noindent{\small {\small\bf Theorem 4.4.}
Let $u_{+}<u_{-}.$    For any fixed $\gamma \in(0, 1)$, assume that $(\rho_{\gamma}(t,x),m_{\gamma}(t,x))=(\rho_{\gamma}(t,x),\rho_{\gamma}(t,x) u_{\gamma}(t,x))$ is a Riemann
solution containing a shock wave and a contact discontinuity of (1.1) and (1.2) with the  Riemann initial data (1.3).
 Then, as $\gamma\rightarrow0$, $(\rho_{\gamma}(t,x),m_{\gamma}(t,x))$ will converge to
 $$(\rho(t,x),m(t,x))=(\rho_{0}(t,x)+w_{1}(t)\delta_{S},\rho_{0}(t,x)u_{0}(t,x)+w_{2}(t)\delta_{S}),$$
in the sense of distributions, and the singular parts of the limit functions
$\rho(t,x)$ and $m(t,x)$  are  a $\delta$-measure with  weights
 $$ w_{1}(t)=t(\sigma[\rho]-[\rho u])=\rho_{-}(u_{-}-u_{+})t,\,\,  \mathrm{and} \,\,\, w_{2}(t)=t(\sigma[\rho u]-[\rho u^{2}]),$$
 respectively, where $\sigma=u_{+}.$

\vskip 0.1in
\noindent{\small {\small\bf Proof.} (1) Set $\xi=\frac{x}{t}.$  Then for any fixed $\gamma \in(0,1)$, the Riemann solution  containing a shock wave and a contact discontinuity of  (1.1) and  (1.2) can be written as $$(\rho_{\gamma},u_{\gamma})(\xi)=\left\{\begin{array}{ll} (\rho_{-},
u_{-}),\,\,\,\,\xi< \sigma_{1},\\(\rho_{\ast},
u_{\ast}),\,\,\,\,\sigma_{1}< \xi < \sigma_{2},\\(\rho_{+},
u_{+}),\,\,\,\,\xi> \sigma_{2}.\end{array} \right.\eqno{}$$
From (3.2), we have the following weak formulations:
$$\int_{-\infty}^{+\infty}\rho_{\gamma}(\xi)(u_{\gamma}(\xi)-\xi)\varphi'(\xi)d\xi
-\int_{-\infty}^{+\infty}\rho_{\gamma}(\xi)\varphi(\xi)d\xi=0,
\eqno{(4.9)}$$
$$\int_{-\infty}^{+\infty}\rho_{\gamma}(\xi)u_{\gamma}(\xi)(u_{\gamma}(\xi)-\xi)\varphi'(\xi)d\xi
+\int_{-\infty}^{+\infty}(\rho_{\gamma}(\xi))^{\gamma+1}(u_{\gamma}(\xi)-\xi)\varphi'(\xi)d\xi
$$$$
-\int_{-\infty}^{+\infty}\left(\rho_{\gamma}(\xi)u_{\gamma}(\xi)+(\rho_{\gamma}(\xi))^{\gamma+1}\right)\varphi(\xi)d\xi=0,\eqno{(4.10)}
$$
for any $\varphi(\xi)\in C_{0}^{+\infty}(R)$.

(2) For the first integral on the left-hand side of (4.9),  using the method of integration by parts, we can derive
$$\int_{-\infty}^{+\infty}\rho_{\gamma}(\xi)(u_{\gamma}(\xi)-\xi)\varphi'(\xi)d\xi
=\left(\int_{-\infty}^{\sigma_{1}}+\int_{\sigma_{2}}^{+\infty}+\int_{\sigma_{1}}^{\sigma_{2}}\right)
\rho_{\gamma}(\xi)(u_{\gamma}(\xi)-\xi)\varphi'(\xi)d\xi
$$
$$=\rho_{-}u_{-}\varphi(\sigma_{1})-\rho_{+}u_{+}\varphi(\sigma_{2})
-\rho_{-}\sigma_{1}\varphi(\sigma_{1})+\rho_{+}\sigma_{2}\varphi(\sigma_{2})+\int_{-\infty}^{\sigma_{1}}\rho_{-}\varphi(\xi)d\xi
$$
$$+\int_{\sigma_{2}}^{+\infty}\rho_{+}\varphi(\xi)d\xi
+\int_{\sigma_{1}}^{\sigma_{2}}\rho_{\ast}(u_{\ast}-\xi)\varphi'(\xi)d\xi
$$
Meanwhile, we have
$$ \int_{\sigma_{1}}^{\sigma_{2}}\rho_{\ast}(u_{\ast}-\xi)\varphi'(\xi)d\xi
=\rho_{\ast}u_{\ast}(\varphi(\sigma_{2})-\varphi(\sigma_{1}))-\rho_{\ast}(\sigma_{2}\varphi(\sigma_{2})-\sigma_{1}\varphi(\sigma_{1}))
+\int_{\sigma_{1}}^{\sigma_{2}}\rho_{\ast}\varphi(\xi)d\xi
$$
$$=\rho_{\ast}(\sigma_{2}-\sigma_{1})\left(u_{\ast}\frac{\varphi(\sigma_{2})-\varphi(\sigma_{1})}{\sigma_{2}-\sigma_{1}}
+\frac{\int_{\sigma_{1}}^{\sigma_{2}}\varphi(\xi)d\xi}{\sigma_{2}-\sigma_{1}}
-\frac{\sigma_{2}\varphi(\sigma_{2})-\sigma_{1}\varphi(\sigma_{1})}{\sigma_{2}-\sigma_{1}}\right).
$$
Then, by Lemma 4.2-4,3, we can obtain
$$\lim_{\gamma\rightarrow0}\int_{\sigma_{1}}^{\sigma_{2}}\rho_{\ast}(u_{\ast}-\xi)\varphi'(\xi)d\xi=0.
$$
 Hence taking the limit $ \gamma\rightarrow0$ in (4.9) leads to
 $$\lim_{\gamma\rightarrow0}\int_{-\infty}^{+\infty}(\rho_{\gamma}(\xi)-\rho_{0}(\xi))\varphi(\xi)d\xi
=(\sigma[\rho]-[\rho u])\varphi(\sigma),\eqno (4.11)$$
where $ (\rho_{0}(\xi),u_{0}(\xi))=(\rho_{\pm},u_{\pm}),~\pm(\xi-\sigma)>0.$

(3) Similarly, we can obtain for (4.10) that
$$
\int_{-\infty}^{+\infty}\rho_{\gamma}(\xi)u_{\gamma}(\xi)(u_{\gamma}(\xi)-\xi)\varphi'(\xi)d\xi
$$
$$=\left(\sigma[\rho u]-[\rho u^{2}]\right)\varphi(\sigma)+\int_{-\infty}^{+\infty}\rho_{0}(\xi)u_{0}(\xi)\varphi(\xi)d\xi
$$
and
$$\int_{-\infty}^{+\infty}(\rho_{\gamma}(\xi))^{\gamma+1}(u_{\gamma}(\xi)-\xi)\varphi'(\xi)d\xi
=\left(\int_{-\infty}^{\sigma_{1}}+\int_{\sigma_{2}}^{+\infty}+\int_{\sigma_{1}}^{\sigma_{2}}\right)(\rho_{\gamma}(\xi))^{\gamma+1}(u_{\gamma}(\xi)-\xi)\varphi'(\xi)d\xi
$$
$$=\rho^{\gamma+1}_{-}u_{-}\varphi(\sigma_{1})-\rho^{\gamma+1}_{+}u_{+}\varphi(\sigma_{2})
-\rho^{\gamma+1}_{-}\sigma_{1}\varphi(\sigma_{1})+\rho^{\gamma+1}_{+}\sigma_{2}\varphi(\sigma_{2})
+\int_{-\infty}^{\sigma_{1}}\rho^{\gamma+1}_{-}\varphi(\xi)d\xi
$$
$$ +\int_{\sigma_{2}}^{+\infty}\rho^{\gamma+1}_{+}\varphi(\xi)d\xi
+\rho^{\gamma}_{\ast}\rho_{\ast}(\sigma_{2}-\sigma_{1})\left(u_{\ast}\frac{\varphi(\sigma_{2})-\varphi(\sigma_{1})}{\sigma_{2}-\sigma_{1}}
-\frac{\sigma_{2}\varphi(\sigma_{2})-\sigma_{1}\varphi(\sigma_{1})}{\sigma_{2}-\sigma_{1}}
+\frac{\int^{\sigma_{2}}_{\sigma_{1}}\varphi(\xi)d\xi}{\sigma_{2}-\sigma_{1}}\right),
$$
which converges to
$$(\sigma[\rho]-[\rho u])\varphi(\sigma)+\int_{-\infty}^{+\infty}\rho_{0}(\xi)\varphi(\xi)d\xi
$$
by Lemma 4.1-4.3.

Thus, following (4.11), we can get
$$ \lim\limits_{\gamma\rightarrow0}\int_{-\infty}^{+\infty}(\rho_{\gamma}(\xi)u_{\gamma}(\xi)-\rho_{0}(\xi)u_{0}(\xi))\varphi(\xi)d\xi
=\left(\sigma[\rho u]-[\rho u^{2}]\right)\varphi(\sigma).\eqno (4.12)
$$

(4) Finally, we study the limits of $\rho_{\gamma}(t,x) $ and $\rho_{\gamma}(t,x)u_{\gamma}(t,x) $ depending on $t$ as $\gamma\rightarrow0$.
Regarding $t$ as a parameter, we can  get from   (4.11) that
$$
\lim_{\gamma\rightarrow0}\int_{-\infty}^{+\infty}(\rho_{\gamma}(\xi)-\rho_{0}(\xi))\varphi(t,\xi t)d\xi
=\lim_{\gamma\rightarrow0}\int_{-\infty}^{+\infty}(\rho_{\gamma}(x/t)-\rho_{0}(x/t))\varphi(t,x)d(x/t)
$$
$$=\frac{1}{t}\lim_{\gamma\rightarrow0}\int_{-\infty}^{+\infty}(\rho_{\gamma}(t,x)-\rho_{0}(t,x))\varphi(t,x)dx
=(\sigma[\rho]-[\rho u])\varphi(t, \sigma t).\eqno (4.13)
$$
Then multiplying (4.13) by $t$ and taking integration, we have
$$\lim_{\gamma\rightarrow0}\int_{0}^{+\infty}\int_{-\infty}^{+\infty}(\rho_{\gamma}(t,x)-\rho_{0}(t,x))\varphi(t,x)dx dt
=\int_{0}^{+\infty}t(\sigma[\rho]-[\rho u])\varphi(t,\sigma t)dt
$$
in which by definition (2.3), we have
$$\int_{0}^{+\infty}t(\sigma[\rho]-[\rho u])\varphi(t,\sigma t)dt=
\langle w_{1}(\cdot)\delta_{S},\varphi(\cdot,\cdot)\rangle.\eqno (4.14)
$$
where
 $$ w_{1}(t)=t(\sigma[\rho]-[\rho u])=\rho_{-}(u_{-}-u_{+})t.$$

In the same way, we can derive from (4.12) that
$$\lim_{\gamma\rightarrow0}\int_{0}^{+\infty}\int_{-\infty}^{+\infty}(\rho_{\gamma}(t,x)u_{\gamma}(t,x)-\rho_{0}u_{0}(t,x))\varphi(t,x)dxdt
=\langle w_{2}(\cdot)\delta_{S},\varphi(\cdot,\cdot)\rangle.\eqno (4.15)
$$
where
 $$ w_{2}(t)=t(\sigma[\rho u]-[\rho u^{2}]).$$
The proof is completed. $~~\Box$

 \vskip 0.1in

\sec{\Large\bf 5.\quad  Riemann solutions of the PAR model (1.7)}
In this section, we construct the Riemann solutions of the perturbed Aw-Rascle model (1.7) with initial data (1.3).

The system (1.7) has two eigenvalues
 $$ \overline{\lambda}_{1}=u-\sqrt{(\gamma-1)\rho^{\gamma-1}u},\,\,\,\,\,\,\,\overline{\lambda}_{2}=u+\sqrt{(\gamma-1)\rho^{\gamma-1}u},\eqno{(5.1)}
$$
with the corresponding right eigenvectors
$$ \overrightarrow{r}_{1}=(\rho,-\sqrt{(\gamma-1)\rho^{\gamma-1}u})^{T},\,\,\,\,\,\,
\overrightarrow{r}_{2}=(\rho,\sqrt{(\gamma-1)\rho^{\gamma-1}u})^{T},
$$
satisfying $\nabla\overline{\lambda}_{i}\cdot\overrightarrow{r_{i}}\neq0$  $(i=1,2)$ for  $\rho>0$ and $ (\gamma+1)\sqrt{u}\pm\sqrt{(\gamma-1)\rho^{\gamma-1}}\neq0$.
Thus, this system   is strictly
hyperbolic  and both characteristic fields are genuinely nonlinear for $\rho,u> 0$  and  $1<\gamma<1+\gamma_{2}$  where $\gamma_{2}>0$  is sufficiently small,
which means  the associated waves are either shock waves or rarefaction waves.

 Seeking the self-similar solution $$(\rho,u)(t,x)=(\rho,u)(\xi),\,\,\,\,\xi=\frac{x}{t},$$
 the Riemann problem    (1.7)   and (1.3) is reduced to the following boundary value problem of the ordinary differential equations:
$$ \left\{\begin{array}{ll}
   -\xi\rho_{\xi}+(\rho u)_{\xi}=0,\\
   -\xi\left(\rho u+\frac{\rho^{\gamma}}{\gamma}\right)_{\xi}+(\rho u^{2}+\rho^{\gamma}u)_{\xi}=0,\end{array}\right .\eqno{(5.2)}
$$
with $(\rho,u)(\pm\infty)=(\rho_{\pm},u_{\pm}).$\\
\indent
For any smooth solution, system (5.2) can be written as
$$
\left(\begin{array}{cccc}u-\xi &\rho\\(u-\xi)u+(\gamma u-\xi)\rho^{\gamma-1} &-\xi \rho+2\rho u+\rho^{\gamma}
 \end{array}\right)\left(\begin{array}{cccc}\rho_{\xi}\\u_{\xi}
 \end{array}\right)=0.\eqno{(5.3)}$$
Besides the constant solution $$(\rho,u)(\xi)={\mathrm constant }   \,\,\,\,\,(\rho> 0),$$  it provides  the 1-rarefaction wave
$$R_{1}(\rho_{-},u_{-}):
\left\{
  \begin{array}{ll}
    \xi=\overline{\lambda}_{1}=u-\sqrt{(\gamma-1)\rho^{\gamma-1}u}, \\
    \sqrt{u}-\sqrt{u_{-}}=-\sqrt{\frac{1}{\gamma-1}\rho^{\gamma-1}}+\sqrt{\frac{1}{\gamma-1}\rho_{-}^{\gamma-1}}, \\
    \rho<\rho_{-},u>u_{-},
  \end{array}
\right.
\eqno{(5.4)}  $$
or the  2-rarefaction wave
$$R_{2}(\rho_{-},u_{-}):
\left\{
  \begin{array}{ll}
    \xi=\overline{\lambda}_{2}=u+\sqrt{(\gamma-1)\rho^{\gamma-1}u}, \\
    \sqrt{u}-\sqrt{u_{-}}=\sqrt{\frac{1}{\gamma-1}\rho^{\gamma-1}}-\sqrt{\frac{1}{\gamma-1}\rho_{-}^{\gamma-1}}, \\
    \rho>\rho_{-},u>u_{-}.
  \end{array}
\right.
\eqno{(5.5)}  $$

Differentiating
 the second equation of  (5.4) with respect to $\rho$ yields $$u_{\rho} = -\sqrt{(\gamma-1)\rho^{\gamma-3}u}<0,$$
 and
$$u_{\rho\rho} =\frac{1}{2}\sqrt{\gamma-1}\left(\sqrt{\gamma-1}\rho^{\gamma-3}-(\gamma-3)\sqrt{\rho^{\gamma-5} u}\right)>0,$$
where $\gamma\in(1,3)$,
which mean that for $1<\gamma<3$, the rarefaction wave curve $R_{1}(\rho_{-},u_{-})$ is monotonic decreasing and convex  in the $(\rho, u)$ phase plane
$(\rho, u> 0)$.

Moreover, by differentiating $\rho$ and $u$ with respect to $\xi$ in the first equation of (5.4) and combining $$ u_{\rho} =\frac{u_{\xi}}{\rho_{\xi}}= -\sqrt{(\gamma-1)\rho^{\gamma-3}u},$$ we have
$$1=\left(\frac{\gamma+1}{2}-\frac{\sqrt{(\gamma-1)\rho^{\gamma-1}}}{2\sqrt{u}}\right)u_{\xi}.\eqno{(5.6)}$$
Hence, as $\gamma\in(1,1+\gamma_{0})$ for $\gamma_{0}$ sufficiently small, we have $u_{\xi}>0$, i.e., the set $(\rho,u)$ which can be joined to
$(\rho_{-},u_{-})$ by 1-rarefaction wave is made up of the half-branch of $R_{1}(\rho_{-},u_{-})$ with $u \geq u_{-}$.

With the same way to compute $R_{2}(\rho_{-},u_{-})$, we can gain $u_{\rho}>0,~~u_{\rho\rho}<0,$ and $u_{\xi}>0$, which  means that it is monotonic creasing
and concave for  $1<\gamma<3$ in the $(\rho, u)$ phase plane $(\rho ,u> 0)$ and the set $(\rho,u)$ which can be joined to
$(\rho_{-},u_{-})$ by 2-rarefaction wave is made up of the half-branch of $R_{1}(\rho_{-},u_{-})$ with $u \geq u_{-}$.

 Performing the limit $\rho\rightarrow 0$ in the second equation in (5.4) yields
$$\lim_{\rho\rightarrow 0}\sqrt{u}=\sqrt{u_{-}}-\lim_{\rho\rightarrow 0}\sqrt{\frac{\rho^{\gamma-1}}{\gamma-1}}
+\sqrt{\frac{\rho_{-}^{\gamma-1}}{\gamma-1}}=\sqrt{u_{-}}+\sqrt{\frac{\rho_{-}^{\gamma-1}}{\gamma-1}}.
$$
Then we have
$$\lim_{\rho\rightarrow 0}u=\left(\sqrt{u_{-}}+\sqrt{\frac{\rho_{-}^{\gamma-1}}{\gamma-1}}\right)^{2}=:u_{0}^{\gamma}.\eqno{(5.7)}
$$Thus we conclude that there exists $u_{0}^{\gamma}$
 such that the 1-rarefaction wave
curve $R_{1}(\rho_{-},u_{-})$ intersects   the $u$-axis at the point $(0, u_{0}^{\gamma}).$

 Performing the limit $\rho\rightarrow +\infty$ of the second equation in (5.5) yields
$$ \lim_{\rho\rightarrow +\infty}\sqrt{u}=\sqrt{u_{-}}+\lim_{\rho\rightarrow +\infty}\left(\sqrt{\frac{\rho^{\gamma-1}}{\gamma-1}}
-\sqrt{\frac{\rho_{-}^{\gamma-1}}{\gamma-1}}\right)=+\infty,\eqno{(5.8)}
$$
which implies that $\lim\limits_{\rho\rightarrow +\infty}u=+\infty$

For a bounded discontinuity at $\xi=\overline{\sigma},$ the Rankine-Hugoniot relation
$$\left\{
    \begin{array}{ll}
      -\overline{\sigma}[\rho]+[\rho u]=0, \\
     -\overline{\sigma}[\rho u+\frac{1}{\gamma}\rho^{\gamma}]+[\rho u^{2}+u\rho^{\gamma}]=0,
    \end{array}
  \right.\eqno (5.9)
$$
holds, where $[\rho]=\rho -\rho_{-},$   etc.   Eliminating $\sigma$ from (5.9), we obtain
$$[\rho][\rho u^{2}] - ([\rho u])^{2}=-[\rho][u\rho^{\gamma}] + [\rho u][\frac{1}{\gamma}\rho^{\gamma}].
\eqno {(5.10)}$$
Simplifying (5.10) yields$$(u - u_{-})^{2}= (\frac{1}{\rho_{-}}- \frac{1}{\rho})(u_{}\rho^{\gamma}-u_{-}\rho_{-}^{\gamma})
-\frac{1}{\gamma\rho\rho_{-}}(\rho u-\rho_{-}u_{-})(\rho^{\gamma}-\rho_{-}^{\gamma}),\eqno {(5.11)}$$
i.e.,$$(u - u_{-})^{2}= \frac{\rho-\rho_{-}}{\rho\rho_{-}}u_{-}(\rho^{\gamma}-\rho_{-}^{\gamma})+\frac{\rho-\rho_{-}}{\rho\rho_{-}}\rho^{\gamma}(u-u_{-})
-\frac{1}{\gamma\rho\rho_{-}}u_{-}(\rho-\rho_{-})(\rho^{\gamma}-\rho_{-}^{\gamma})
-\frac{1}{\gamma\rho\rho_{-}}\rho(u-u_{-})(\rho^{\gamma}-\rho_{-}^{\gamma}).$$
Therefore,$$(\frac{u - u_{-}}{\rho-\rho_{-}})^{2}=(1-\frac{1}{\gamma}) \frac{u_{-}}{\rho\rho_{-}}(\frac{\rho^{\gamma}-\rho_{-}^{\gamma}}{\rho-\rho_{-}})+\frac{u-u_{-}}{\rho-\rho_{-}}\bigg(\frac{\rho^{\gamma-1}}{\rho_{-}}
-\frac{1}{\gamma\rho_{-}}(\frac{\rho^{\gamma}-\rho_{-}^{\gamma}}{\rho-\rho_{-}})\bigg).\eqno {(5.12)}$$
Set $\alpha = \frac{u - u_{-}}{\rho-\rho_{-}}$. Then (5.12)  can be simplified as
$$\alpha^{2}-\bigg(\frac{\rho^{\gamma-1}}{\rho_{-}}
-\frac{1}{\gamma\rho_{-}}(\frac{\rho^{\gamma}-\rho_{-}^{\gamma}}{\rho-\rho_{-}})\bigg)\alpha-(1-\frac{1}{\gamma}) \frac{u_{-}}{\rho\rho_{-}}(\frac{\rho^{\gamma}-\rho_{-}^{\gamma}}{\rho-\rho_{-}})=0.\eqno {}$$This is a quadratic form in $\alpha$ and we can solve this to obtain
      $$ \frac{u - u_{-}}{\rho-\rho_{-}}=\frac{1}{2\rho\rho_{-}}
\bigg(\rho^{\gamma}
-\frac{\rho}{\gamma}(\frac{\rho^{\gamma}-\rho_{-}^{\gamma}}{\rho-\rho_{-}})\bigg)\pm\sqrt{\frac{1}{4\rho^{2}\rho_{-}^{2}}
\bigg(\rho^{\gamma}
-\frac{\rho}{\gamma}(\frac{\rho^{\gamma}-\rho_{-}^{\gamma}}{\rho-\rho_{-}})\bigg)^{2}+(1-\frac{1}{\gamma}) \frac{u_{-}}{\rho\rho_{-}}(\frac{\rho^{\gamma}-\rho_{-}^{\gamma}}{\rho-\rho_{-}})},
     \eqno (5.13)$$
where  $(\rho_{-},u_{-})$ and $(\rho_{},u_{})$ are the shock speed, the left state and the right state, respectively.

$\textbf{1-shock~wave}~S_{1}(\rho_{-},u_{-})$\textbf{:}

The classical Lax
 entropy conditions   imply  that the propagation speed $\overline{\sigma}$ for the 1-shock wave
has to be satisfied with$$\overline{\lambda}_{1}(\rho,u)<\overline{\sigma}<\overline{\lambda}_{1}(\rho_{-},u_{-}).\eqno (5.14)
$$
From the first equation of (5.9), we have
$$\overline{\sigma}=\frac{\rho u -\rho_{-}u_{-}}{\rho-\rho_{-}} = u_{-}+\frac{\rho}{\rho-\rho_{-}}(u-u_{-}).$$
Then, it follows from the right inequality of (5.14) that
$$\frac{\rho}{\rho-\rho_{-}}(u-u_{-})<-\sqrt{(\gamma-1)\rho_{-}^{\gamma-1}u_{-}}<0,\eqno (5.15)
$$
which implies that $u-u_{-}$ and $\rho-\rho_{-}$ have different signs.  Similarly, for the left inequality of (5.14), we can gain
 $$\frac{\rho_{-}}{\rho-\rho_{-}}(u-u_{-})>-\sqrt{(\gamma-1)\rho^{\gamma-1}u}.\eqno (5.16)
$$
 Combining (5.15) and (5.16), it is easy to get
$$-\sqrt{(\gamma-1)\rho^{\gamma+1}u}<\frac{\rho\rho_{-}}{\rho-\rho_{-}}(u-u_{-})<-\sqrt{(\gamma-1)\rho_{-}^{\gamma+1}u_{-}},
$$
which indicates that  $\rho > \rho_{-}, u_{-}> u$, and the minus sign is taken in (5.13) for 1-shock wave.  Hence given a left state $(\rho_{-},u_{-}),$
 the 1-shock wave curve
$S_{1}(\rho_{-},u_{-}) $  in the phase plane which is the set of states that can be connected on the right by a 1-shock
is as follows
$$ u - u_{-}=\frac{1}{2\rho\rho_{-}}
\bigg(\rho^{\gamma}(\rho-\rho_{-})
-\frac{\rho}{\gamma}(\rho^{\gamma}-\rho_{-}^{\gamma})\bigg)-\sqrt{\frac{1}{4\rho^{2}\rho_{-}^{2}}
\bigg(\rho^{\gamma}({\rho-\rho_{-}})
-\frac{\rho}{\gamma}(\rho^{\gamma}-\rho_{-}^{\gamma})\bigg)^{2}+(1-\frac{1}{\gamma}) \frac{u_{-}}{\rho\rho_{-}}(\rho^{\gamma}-\rho_{-}^{\gamma})(\rho-\rho_{-})}, ~~\rho>\rho_{-}.
     \eqno (5.17)$$

$\textbf{2-shock~wave~}S_{2}(\rho_{-},u_{-}):$

The propagation speed $\overline{\sigma}$ for the $2$-shock wave should satisfy
$$\overline{\lambda}_{2}(\rho,u)<\overline{\sigma}<\overline{\lambda}_{2}(\rho_{-},u_{-}).\eqno (5.18)
$$
With the similar calculations to the 1-shock wave, we have the the 2-shock curve $S_{2}(\rho_{-},u_{-}): $
$$ u - u_{-}=\frac{1}{2\rho\rho_{-}}
\bigg(\rho^{\gamma}(\rho-\rho_{-})
-\frac{\rho}{\gamma}(\rho^{\gamma}-\rho_{-}^{\gamma})\bigg)-\sqrt{\frac{1}{4\rho^{2}\rho_{-}^{2}}
\bigg(\rho^{\gamma}({\rho-\rho_{-}})
-\frac{\rho}{\gamma}(\rho^{\gamma}-\rho_{-}^{\gamma})\bigg)^{2}+(1-\frac{1}{\gamma}) \frac{u_{-}}{\rho\rho_{-}}(\rho^{\gamma}-\rho_{-}^{\gamma})(\rho-\rho_{-})}, ~~\rho<\rho_{-}.
     \eqno (5.19)$$

Differentiating $u$ with respect to $\rho$ in the second equation of (5.11) gives
 that for $\rho>\rho_{-}$,
$$\rho\rho_{-}I_{1}u_{\rho}=I_{2},\eqno{(5.20)}$$
where
$$I_{1}=2(u-u_{-})-\frac{\gamma-1}{\gamma}\left(\frac{1}{\rho_{-}}-\frac{1}{\rho}\right)\rho^{\gamma}
+\frac{1}{\gamma}\left(\rho^{\gamma-1}-\rho_{-}^{\gamma-1}\right)<0,
$$
$$I_{2}=\frac{\gamma-1}{\gamma\rho}\left((\rho^{\gamma}-\rho_{-}^{\gamma})\rho_{-}u_{-}+\gamma(\rho-\rho_{-})\rho^{\gamma}u\right),
$$
which gives $u_{\rho}<0$ for $\gamma\in(1,1+\gamma_{0})$ where $\gamma_{0}$ sufficiently small, which indicates that  the 1-shock wave curve $S_{1}(\rho_{-},u_{-})$ is monotonic decreasing in the region $\rho > \rho_{-}$ in the $(\rho,u)$ phase plane. Moreover, letting $u=0$ in (5.11), it is easy to get
$$u_{-}=\sqrt{\frac{1}{\gamma}\left(u_{-}(\rho^{\gamma-1}-\rho_{-}^{\gamma-1})-(\gamma-1)\left(\frac{1}{\rho_{-}}-\frac{1}{\rho}\right)
\rho_{-}^{\gamma}u_{-}\right)}.\eqno{(5.21)}
$$
Setting
$$ f(\rho)=u_{-}-\sqrt{\frac{1}{\gamma}\left(u_{-}(\rho^{\gamma-1}-\rho_{-}^{\gamma-1})-(\gamma-1)\left(\frac{1}{\rho_{-}}-\frac{1}{\rho}\right)
\rho_{-}^{\gamma}u_{-}\right)}.
$$
Then $f(\rho_{-})f(+\infty)<0$, and $f (\rho)$ is continuous with respect to $\rho$. Therefore, there exists $\rho_{0}\in(\rho_{-},+\infty)$ such that
$f(\rho_{0})=0$, which implies that  the 1-shock wave curve $S_{1}(\rho_{-},u_{-})$ intersects
with the $\rho$-axis at a point.

Similarly,  we can
get $u_{\rho}>0$ for the 2-shock wave for  for $\gamma\in(1,1+\gamma_{0})$ where $\gamma_{0}$ sufficiently small, which indicates that  the 2-shock wave curve $S_{2}(\rho_{-},u_{-})$ is monotonic increasing in the region $\rho < \rho_{-}$ in the $(\rho,u)$ phase plane.
 From (5.19), it is not difficult to check that
that $\lim\limits_{\rho\rightarrow 0^{+}}u= -\infty$ for the 2-shock wave curve $S_{2}^{\gamma} (\rho_{-},u_{-})$, which implies that curve $S_{2}^{\gamma} (\rho_{-},u_{-})$
has the $u$-axis  as its asymptotic line.

In the $(\rho, u) $ phase plane $( \rho, u \geq¡Ý 0)$, through a given point $(\rho_{-}, u_{-})$, we draw the elementary wave curves.  We find that the elementary wave curves divide the quarter phase plane $( \rho, u \geq¡Ý 0)$ into five regions,
 see Fig. 2.  According to the right
state $(\rho_{+},u_{+})$ in the different regions, one can construct the unique global solution to the Riemann problem (1.7)  and (1.3) as follows:

(1) $(\rho_{+},u_{+})\in I(\rho_{-},u_{-}):$ $R_{1}+R_{2};$

 (2)$(\rho_{+},u_{+})\in II(\rho_{-},u_{-}):$ $S_{1}+R_{2};$

 (3)$(\rho_{+},u_{+})\in III(\rho_{-},u_{-}):$ $R_{1}+S_{2};$

 (4)$(\rho_{+},u_{+})\in IV(\rho_{-},u_{-}):$ $S_{1}+S_{2};$

  (5)$(\rho_{+},u_{+})\in V(\rho_{-},u_{-}):$ $R_{1}+\mathrm{Vac}+R_{2},$\\
 where   ``+" means ``followed by".

\hspace{65mm}\setlength{\unitlength}{0.8mm}\begin{picture}(80,66)
\put(-32,-0.2){\vector(0,2){50}}\put(1.7,9.6){$\cdot$}
\put(28.5,-5){$u_{0}^{\gamma}$
}
 \put(-50,0){\vector(2,0){115}}  \put(-35,49){$\rho$}\put(-35,-4){$0$}
\put(66,-1){$u$}
\put(23,44){$R_{2}$}
\put(-25,3){$S_{2}$}\put(55,40){$R_{2}$}\put(20,3){$R_{1}$}

\put(1,13){$(\rho_{-},
u_{-})$}
\put(44,6){$$ }\put(2,32){}\qbezier(-32,1)(17.5,4)(32,52)\qbezier(30,0)(-14,12)(-32,42)\qbezier(30,0)(47.5,4)(62,42)
\put(-25,32){$S_{1}$}\put(-10,3.5){$III (\rho_{-},u_{-})$}
\put(-30,10){$IV(\rho_{-},u_{-})$}\put(50,8){$V(\rho_{-},u_{-})$}
\put(26,18){$I (\rho_{-},u_{-})$}\put(-10,35){$II (\rho_{-},u_{-})$}
\end{picture}}
\vspace{0.6mm}  \vskip 0.2in \centerline{\bf Fig. 2.\,   Curves of  elementary waves.
   } \vskip 0.1in \indent

\baselineskip 15pt
 \sec{\Large\bf 6.\quad   Limits of Riemann solutions of (1.7)}
In this section, we  study the limiting behavior of the Riemann solutions of (1.7)   as $\gamma$ goes to one, that is, the formation of delta shock and the vacuum states as $\gamma\rightarrow 1$, respectively in the case $u_{-}>u_{+}$ and in the case $u_{-}<u_{+}$.

\baselineskip 15pt
 \sec{\Large\bf 6.1.\quad   Formation of delta  shock wave }
In this subsection, we study the formation of $\delta$-shock in the Riemann problem (1.7) and (1.3) when
$u_{-}>u_{+}$ as $\gamma\rightarrow 1$.

 \vskip 0.1in
\noindent{\small {\small\bf Lemma 6.1.}
If $u_{+}<u_{-}$, then there is a sufficiently small $\gamma_{0}>0$ such that $(\rho_{+},u_{+})\in IV(\rho_{-},u_{-})$ as $1<\gamma <1 +\gamma_{0}$.

\vskip 0.1in
\noindent{\small {\small\bf Proof.} If $\rho_{+} = \rho_{-}$, then $(\rho_{+}, u_{+})\in IV( \rho_{-}, u_{-})$  for any $\gamma\in (1,3)$. Thus, we only need to consider the
case $\rho_{+} \neq\rho_{-}$.

It can be derived from (5.17) and (5.19) that all possible states $(\rho, u)$ that can be connected to the
left state $(\rho_{-}, u_{-})$ on the right by a 1-shock wave $S_{1}$ or a 2-shock wave $S_{2}$ should satisfy
$$ S_{1}:  u = u_{-}+\frac{1}{2\rho\rho_{-}}
\bigg(\rho^{\gamma}(\rho-\rho_{-})
-\frac{\rho}{\gamma}(\rho^{\gamma}-\rho_{-}^{\gamma})\bigg)$$$$-(\rho-\rho_{-})\sqrt{\frac{1}{4\rho^{2}\rho_{-}^{2}}
\bigg(\rho^{\gamma}
-\frac{\rho}{\gamma}\bigg(\frac{\rho^{\gamma}-\rho_{-}^{\gamma}}{\rho-\rho_{-}}\bigg)\bigg)^{2}+(1-\frac{1}{\gamma}) \frac{u_{-}}{\rho\rho_{-}}\bigg(\frac{\rho^{\gamma}-\rho_{-}^{\gamma}}{\rho-\rho_{-}}\bigg)}, ~~\rho>\rho_{-},
     \eqno (6.1)$$$$ S_{2}:  u = u_{-}+\frac{1}{2\rho\rho_{-}}
\bigg(\rho^{\gamma}(\rho-\rho_{-})
-\frac{\rho}{\gamma}(\rho^{\gamma}-\rho_{-}^{\gamma})\bigg)$$$$+(\rho-\rho_{-})\sqrt{\frac{1}{4\rho^{2}\rho_{-}^{2}}
\bigg(\rho^{\gamma}
-\frac{\rho}{\gamma}\bigg(\frac{\rho^{\gamma}-\rho_{-}^{\gamma}}{\rho-\rho_{-}}\bigg)\bigg)^{2}+(1-\frac{1}{\gamma}) \frac{u_{-}}{\rho\rho_{-}}\bigg(\frac{\rho^{\gamma}-\rho_{-}^{\gamma}}{\rho-\rho_{-}}\bigg)}, ~~\rho<\rho_{-}.
     \eqno (6.2)$$

If $\rho_{+} \neq\rho_{-}$ and $(\rho_{+}, u_{+})\in IV( \rho_{-}, u_{-})$, then from Fig. 1, (6.1) and (6.2), we have$$  u_{+}< u_{-}+\frac{1}{2\rho_{+}\rho_{-}}
\bigg(\rho_{+}^{\gamma}(\rho_{+}-\rho_{-})
-\frac{\rho_{+}}{\gamma}(\rho_{+}^{\gamma}-\rho_{-}^{\gamma})\bigg)$$$$-(\rho_{+}-\rho_{-})\sqrt{\frac{1}{4\rho_{+}^{2}\rho_{-}^{2}}
\bigg(\rho_{+}^{\gamma}
-\frac{\rho_{+}}{\gamma}\bigg(\frac{\rho_{+}^{\gamma}-\rho_{-}^{\gamma}}{\rho_{+}-\rho_{-}}\bigg)\bigg)^{2}+(1-\frac{1}{\gamma}) \frac{u_{-}}{\rho_{+}\rho_{-}}\bigg(\frac{\rho_{+}^{\gamma}-\rho_{-}^{\gamma}}{\rho_{+}-\rho_{-}}\bigg)}, ~~\rho_{+}>\rho_{-},
     \eqno (6.3)$$$$   u_{+} < u_{-}+\frac{1}{2\rho_{+}\rho_{-}}
\bigg(\rho_{+}^{\gamma}(\rho_{+}-\rho_{-})
-\frac{\rho_{+}}{\gamma}(\rho_{+}^{\gamma}-\rho_{-}^{\gamma})\bigg)$$$$+(\rho_{+}-\rho_{-})\sqrt{\frac{1}{4\rho_{+}^{2}\rho_{-}^{2}}
\bigg(\rho_{+}^{\gamma}
-\frac{\rho_{+}}{\gamma}\bigg(\frac{\rho_{+}^{\gamma}-\rho_{-}^{\gamma}}{\rho_{+}-\rho_{-}}\bigg)\bigg)^{2}+(1-\frac{1}{\gamma}) \frac{u_{-}}{\rho_{+}\rho_{-}}\bigg(\frac{\rho_{+}^{\gamma}-\rho_{-}^{\gamma}}{\rho_{+}-\rho_{-}}\bigg)}, ~~\rho_{+}<\rho_{-},
     \eqno (6.4)$$
which implies that
$$ \sqrt{\frac{1}{4\rho_{+}^{2}\rho_{-}^{2}}
\bigg(\rho_{+}^{\gamma}
-\frac{\rho_{+}}{\gamma}\bigg(\frac{\rho_{+}^{\gamma}-\rho_{-}^{\gamma}}{\rho_{+}-\rho_{-}}\bigg)\bigg)^{2}+(1-\frac{1}{\gamma}) \frac{u_{-}}{\rho_{+}\rho_{-}}\bigg(\frac{\rho_{+}^{\gamma}-\rho_{-}^{\gamma}}{\rho_{+}-\rho_{-}}\bigg)}$$
$$-\frac{1}{2}\bigg|\frac{1}{\rho_{-}}-\frac{1}{\rho_{+}}\bigg|\bigg(\frac{\rho_{+}^{\gamma}}{\rho_{+}-\rho_{-}}
-\frac{\rho_{+}(\rho_{+}^{\gamma}-\rho_{-}^{\gamma})}{\gamma(\rho_{+}-\rho_{-})^{2}}\bigg)< \frac{u_{-}-u_{+}}{|\rho_{+}-\rho_{-}|}. \eqno{(6.5)}$$
Since
$$  \lim_{{\gamma\rightarrow1}}\Bigg( \sqrt{\frac{1}{4\rho_{+}^{2}\rho_{-}^{2}}
\bigg(\rho_{+}^{\gamma}
-\frac{\rho_{+}}{\gamma}\bigg(\frac{\rho_{+}^{\gamma}-\rho_{-}^{\gamma}}{\rho_{+}-\rho_{-}}\bigg)\bigg)^{2}+(1-\frac{1}{\gamma}) \frac{u_{-}}{\rho_{+}\rho_{-}}\bigg(\frac{\rho_{+}^{\gamma}-\rho_{-}^{\gamma}}{\rho_{+}-\rho_{-}}\bigg)}$$
$$-\frac{1}{2}\bigg|\frac{1}{\rho_{-}}-\frac{1}{\rho_{+}}\bigg|\bigg(\frac{\rho_{+}^{\gamma}}{\rho_{+}-\rho_{-}}
-\frac{\rho_{+}(\rho_{+}^{\gamma}-\rho_{-}^{\gamma})}{\gamma(\rho_{+}-\rho_{-})^{2}}\bigg)\Bigg)=0,
  \eqno{(6.6)}$$
it follows that there exists $\gamma_{0}>0$ small enough such that, when $ 1<\gamma <1 +\gamma_{0}$,  we have
$$ \sqrt{\frac{1}{4\rho_{+}^{2}\rho_{-}^{2}}
\bigg(\rho_{+}^{\gamma}
-\frac{\rho_{+}}{\gamma}\bigg(\frac{\rho_{+}^{\gamma}-\rho_{-}^{\gamma}}{\rho_{+}-\rho_{-}}\bigg)\bigg)^{2}+(1-\frac{1}{\gamma}) \frac{u_{-}}{\rho_{+}\rho_{-}}\bigg(\frac{\rho_{+}^{\gamma}-\rho_{-}^{\gamma}}{\rho_{+}-\rho_{-}}\bigg)}$$
$$-\frac{1}{2}\bigg|\frac{1}{\rho_{-}}-\frac{1}{\rho_{+}}\bigg|\bigg(\frac{\rho_{+}^{\gamma}}{\rho_{+}-\rho_{-}}
-\frac{\rho_{+}(\rho_{+}^{\gamma}-\rho_{-}^{\gamma})}{\gamma(\rho_{+}-\rho_{-})^{2}}\bigg)< \frac{u_{-}-u_{+}}{|\rho_{+}-\rho_{-}|}. \eqno{}$$
  Then,  it is obvious
that $(\rho_{+}, u_{+})\in IV( \rho_{-}, u_{-})$ when $ 1<\gamma <1 +\gamma_{0}$.
 The proof is completed. $~\Box$

\vskip 0.1in
According to the relation  (5.11), for a given state $(\rho_{-},u_{-})$,  the shock curves  $S_{1}(\rho_{-},u_{-})$ and  $S_{2} (\rho_{-},u_{-})$  can also be expressed as below:
$$ u-u_{-}=-\sqrt{\frac{1}{\gamma}\left((\gamma-1)\left(\frac{1}{\rho_{-}}-\frac{1}{\rho}\right)
           (\rho^{\gamma}u-\rho_{-}^{\gamma}u_{-})+(u_{-}-u)(\rho^{\gamma-1}-\rho_{-}^{\gamma-1})\right)}, ~~u<u_{-},
\eqno{(6.7)}  $$
with $\rho>\rho_{-}$ for a 1-shock curve $S_{1} (\rho_{-},u_{-})$ , and $\rho<\rho_{-} $ for a 2-shock curve  $S_{2} (\rho_{-},u_{-})$.

When  $1<\gamma <1+\gamma_{0}$, namely $(\rho_{+}, u_{+})\in IV( \rho_{-}, u_{-})$, suppose that $(\rho_{\ast},
u_{\ast})$ is the intermediate state
connected with $(\rho_{-}, u_{-})$ by a 1-shock wave $S_{1}$ with the speed $\sigma_{1}$, and $(\rho_{+}, u_{+})$ by a 2-shock wave $S_{2}$  with the speed $\sigma_{2},$
then it follows from
 (6.7) that$$ u_{\ast}-u_{-}=-\sqrt{\frac{1}{\gamma}\left((\gamma-1)\left(\frac{1}{\rho_{-}}-\frac{1}{\rho_{\ast}}\right)
(\rho_{\ast}^{\gamma}u_{\ast}-\rho_{-}^{\gamma}u_{-})+(u_{-}-u_{\ast})(\rho_{\ast}^{\gamma-1}-\rho_{-}^{\gamma-1})\right)}, ~\rho_{\ast}>\rho_{-},u_{\ast}<u_{-},\eqno{(6.8)}
$$
$$u_{+}-u_{\ast}=-\sqrt{\frac{1}{\gamma}\left((\gamma-1)\left(\frac{1}{\rho_{\ast}}-\frac{1}{\rho_{+}}\right)
(\rho_{+}^{\gamma}u_{+}-\rho_{\ast}^{\gamma}u_{\ast})+(u_{\ast}-u_{+})(\rho_{+}^{\gamma-1}-\rho_{\ast}^{\gamma-1})\right)}, ~\rho_{\ast}>\rho_{+},u_{\ast}>u_{+},\eqno{(6.9)}
$$
with the shock speed
$$\overline{\sigma}_{1}=\frac{\rho_{\ast}u_{\ast}-\rho_{-}u_{-}}{\rho_{\ast}-\rho_{-}},
~\overline{\sigma}_{2}=\frac{\rho_{+}u_{+}-\rho_{\ast}u_{\ast}}{\rho_{+}-\rho_{\ast}},\eqno{(6.10)}
$$
respectively. In this case, the Riemann solution is
$$ (\rho, u)(t, x) =\left\{\begin{array}{ll} (\rho_{-},
u_{-}),\,\,\,\,x< \overline{\sigma}_{1}t,\\(\rho_{\ast},
u_{\ast}),\,\,\,\,\overline{\sigma}_{1}t< x < \overline{\sigma}_{2}t,\\(\rho_{+},
u_{+}),\,\,\,\,x> \overline{\sigma}_{2}t.\end{array} \right.\eqno{(6.11)}$$

 Based on (6.8) and (6.9), we can get that
$$u_{-}-u_{+}=\sqrt{\frac{1}{\gamma}\left((\gamma-1)\left(\frac{1}{\rho_{-}}-\frac{1}{\rho_{\ast}}\right)
(\rho_{\ast}^{\gamma}u_{\ast}-\rho_{-}^{\gamma}u_{-})+(u_{-}-u_{\ast})(\rho_{\ast}^{\gamma-1}-\rho_{-}^{\gamma-1})\right)}
$$
$$+\sqrt{\frac{1}{\gamma}\left((\gamma-1)\left(\frac{1}{\rho_{+}}-\frac{1}{\rho_{\ast}}\right)
(\rho_{\ast}^{\gamma}u_{\ast}-\rho_{+}^{\gamma}u_{+})+(u_{+}-u_{\ast})(\rho_{\ast}^{\gamma-1}-\rho_{+}^{\gamma-1})\right)}, \,\,\,\,   \rho_{\ast}>\rho_{\pm}, u_{+}<u_{\ast}<u_{-}.\eqno{(6.12)}
$$

 Then we have the following lemmas.

 \vskip 0.1in

\vskip 0.1in
\noindent{\small {\small\bf Lemma 6.2.}
$\lim\limits_{\gamma\rightarrow 1}\rho_{\ast}=+\infty,$ and $ \lim\limits_{\gamma\rightarrow 1}(\gamma-1)\rho_{\ast}^{\gamma}u_{\ast}=:a
=\bigg(\frac{\sqrt{\rho_{-}\rho_{+}}}{\sqrt{\rho_{-}}+\sqrt{\rho_{+}}}(u_{-}-u_{+})\bigg)^{2}$.

\vskip 0.1in
\noindent{\small {\small\bf Proof.}
 Let $ \lim\limits_{\gamma\rightarrow1}\inf\rho_{\ast}=\alpha$, and $\lim\limits_{\gamma\rightarrow1}\sup\rho_{\ast}=\beta$.

If $ \alpha<\beta$ , then by the continuity of $\rho_{\ast}(\gamma)$, there exists a sequence $  \{\gamma_{n}\}_{n=1}^{\infty}\subseteq(1,3)$
such that
$$ \lim_{n\rightarrow +\infty}\gamma_{n}=1,\,\,\mathrm{ and}\,\, \lim_{n\rightarrow +\infty}\rho_{\ast}(\gamma_{n})=c,$$
for  some $ c\in(\alpha,\beta).$ Then substituting the sequence  into the right hand side of (6.12),   taking the limit $n\rightarrow +\infty$, and noting
 $u_{+} < u_{\ast}
< u_{-} $  in mind,
 we have
$$\lim_{n\rightarrow+\infty}\frac{1}{\gamma_{n}}\left((\gamma_{n}-1)\left(\frac{1}{\rho_{\pm}}-\frac{1}{\rho_{\ast}(\gamma_{n})}\right)
\left(\left(\rho_{\ast}(\gamma_{n})\right)^{\gamma_{n}}u_{\ast}-\rho_{\pm}^{\gamma_{n}}u_{\pm}\right)+(u_{\pm}-u_{\ast})
\left(\left(\rho_{\ast}(\gamma_{n})\right)^{\gamma_{n}-1}-\rho_{\pm}^{\gamma_{n}-1}\right)\right)=0.
\eqno{(6.13)} $$
Thus,  we can obtain from (6.12) that
$$ u_{-}-u_{+}=0,$$
which contradicts with the assumption  $u_{-}>u_{+}$.
Then we must have $\alpha=\beta$, which  means $\lim\limits_{\gamma\rightarrow1}\rho_{\ast}(\gamma)=\alpha.$

If  $\alpha\in(0,+\infty),$ then  we  can also get a contradiction when taking limit in (6.12). Hence $\alpha=0 $ or $ \alpha=+\infty$. By the condition
$\rho_{\ast}>\max\{\rho_{-},\rho_{+}\}$, it is easy to see that $\lim\limits_{\gamma\rightarrow 1}\rho_{\ast}(\gamma)=\alpha=+\infty.$

Next  taking the limit $ \gamma\rightarrow 1$  in (6.12), we have
$$u_{-}-u_{+}=\sqrt{\lim_{\gamma\rightarrow 1}(\gamma-1)\rho_{\ast}^{\gamma}u_{\ast}(\frac{1}{\rho_{-}}-\frac{1}{\rho_{\ast}})}
+\sqrt{\lim_{\gamma\rightarrow 1}(\gamma-1)\rho_{\ast}^{\gamma}u_{\ast}(\frac{1}{\rho_{+}}-\frac{1}{\rho_{\ast}})}=:(\sqrt{\frac{1}{\rho_{-}}}+\sqrt{\frac{1}{\rho_{+}}})\sqrt{a},$$
from which we  can get
$ a=\bigg(\frac{\sqrt{\rho_{-}\rho_{+}}}{\sqrt{\rho_{-}}+\sqrt{\rho_{+}}}(u_{-}-u_{+})\bigg)^{2}.$
 The proof is completed. $~\Box$

\vskip 0.1in
\noindent{\small {\small\bf Lemma 6.3.}

$$ \lim_{\gamma\rightarrow 1}\overline{\sigma}_{1}=\lim_{\gamma\rightarrow 1}\overline{\sigma}_{2}=\lim_{\gamma\rightarrow1}u_{\ast}=\sigma,\eqno{(6.14)}
$$and
$$ \lim\limits_{\gamma\rightarrow 1}\int_{\overline{\sigma}_{1}}^{\overline{\sigma}_{2}}\rho_{*}d\xi =\sigma[\rho]-[\rho u], \eqno{(6.15)}
$$
where $\sigma=\frac{\sqrt{\rho_{-}}u_{-}+\sqrt{\rho_{+}}u_{+}}{\sqrt{\rho_{-}}+\sqrt{\rho_{+}}}$.

\vskip 0.1in
\noindent{\small {\small\bf Proof.} From (6.8)-(6.10) and  Lemma 6.2, we immediately get
$$\lim_{\gamma\rightarrow 1}u_{\ast}
=u_{-}-\lim_{\gamma\rightarrow 1}\sqrt{\frac{1}{\gamma}\left((\gamma-1)\left(\frac{1}{\rho_{-}}-\frac{1}{\rho_{\ast}}\right)
(\rho_{\ast}^{\gamma}u_{\ast}-\rho_{-}^{\gamma}u_{-})+(u_{-}-u_{\ast})(\rho_{\ast}^{\gamma-1}-\rho_{-}^{\gamma-1})\right)}
$$
$$=u_{-}-\sqrt{\frac{a}{\rho_{-}}}
=u_{-}-\frac{\sqrt{\rho_{-}\rho_{+}}(u_{-}-u_{+})}{\sqrt{\rho_{-}}(\sqrt{\rho_{-}}+\sqrt{\rho_{+}})}=\sigma,
$$
$$\lim_{\gamma\rightarrow1}\overline{\sigma}_{1}=\lim_{\gamma\rightarrow1}\frac{\rho_{\ast}u_{\ast}-\rho_{-}u_{-}}{\rho_{\ast}-\rho_{-}}
=u_{-}+\lim_{\gamma\rightarrow 1}\frac{\rho_{\ast}}{\rho_{-}-\rho_{\ast}}(u_{-}-u{\ast})
=u_{-}-\frac{\sqrt{\rho_{-}\rho_{+}}(u_{-}-u_{+})}{\sqrt{\rho_{-}}(\sqrt{\rho_{-}}+\sqrt{\rho_{+}})}=\sigma,
$$
and
$$\lim_{\gamma\rightarrow1}\overline{\sigma}_{2}=\lim_{\gamma\rightarrow1}\frac{\rho_{+}u_{+}-\rho_{\ast}u_{\ast}}{\rho_{+}-\rho_{\ast}}
=u_{+}+\lim_{\gamma\rightarrow1}\frac{\rho_{\ast}}{\rho_{+}-\rho_{\ast}}(u_{+}-u{\ast})=u_{+}+\sqrt{\frac{a}{\rho_{+}}}
=u_{+}+\frac{\sqrt{\rho_{-}\rho_{+}}(u_{-}-u_{+})}{\sqrt{\rho_{+}}(\sqrt{\rho_{-}}+\sqrt{\rho_{+}})}=\sigma.
$$

From the first equations of the Rankine-Hugoniot relation (5.9) for $S_{1}$ and $S_{2}$,   we have
 $$ \overline{\sigma}_{1}(\rho_{-}-\rho_{\ast})=\rho_{-}u_{-}-\rho_{\ast}u_{\ast},\eqno{(6.16)}$$and $$\overline{\sigma}_{2}(\rho_{\ast}-\rho_{+})=\rho_{\ast}u_{\ast}-\rho_{+}u_{+}. \eqno{(6.17)}$$
By  (6.14), (6.16) and (6.17), we get
 $$ \lim_{\gamma\rightarrow1}\rho_{\ast}(\overline{\sigma}_{2}-\overline{\sigma}_{1})=\lim_{\gamma\rightarrow1}(\rho_{-}u_{-}
-\overline{\sigma}_{1}\rho_{-}+\overline{\sigma}_{2}\rho_{+}-\rho_{+}u_{+})
=\sigma[\rho]-[\rho u],
$$
which implies that$$\lim\limits_{\gamma\rightarrow 1}\int_{\overline{\sigma}_{1}}^{\overline{\sigma}_{2}}\rho_{*}d\xi =\sigma[\rho]-[\rho u].\eqno{}$$
The proof is completed. $~~\Box$

\vskip 0.1in
\noindent{\small {\small\bf Remark 6.1.}  Lemmas 6.2-6.3 show that when $\gamma$ tends to one, the two shock curves $S_{1}$ and $S_{2}$ coincide to form a new delta shock wave, and the delta shock wave speed $\sigma$ is the limit of both the particle velocity $u_{\ast}$ and two shocks' speed $\overline{\sigma}_{1},~\overline{\sigma}_{2}$.
What is more, the intermediate density $\rho_{\ast}$ tend to singular as $\gamma\rightarrow 1$.

\vskip 0.1in
What is more, we will further derive that, when $\gamma\rightarrow 1$, the limit of Riemann solutions of (1.7) with the  Riemann initial data (1.3) under the assumption $u_{+}<u_{-}$  is  a delta shock wave solution of the zero pressure gas dynamics (1.5)  with the same Riemann initial data
$(\rho_{\pm}, u_{\pm})$  in the sense of distributions.

\vskip 0.1in
\noindent{\small {\small\bf Theorem 6.4.}
Let $u_{+}<u_{-}.$  For any fixed $\gamma \in(1, 3)$, assume that $(\rho_{\gamma}(t,x),m_{\gamma}(t,x))=(\rho_{\gamma}(t,x),\rho_{\gamma}(t,x) u_{\gamma}(t,x))$ is a Riemann solution
containing two shocks $S_{1}$ and $S_{2}$ of (1.7)  with the  Riemann initial data (1.3) constructed in Section 5.
 Then, as $\gamma\rightarrow1$, $(\rho_{\gamma}(t,x),m_{\gamma}(t,x))$ will converge to
 $$(\rho(t,x),m(t,x))=(\rho_{0}(t,x)+w_{1}(t)\delta_{S},\rho_{0}(t,x)u_{0}(t,x)+w_{2}(t)\delta_{S}),$$
in the sense of distributions, and the singular parts of the limit functions
$\rho(t,x)$ and $m(t,x)$  are  a $\delta$-measure with  weights
 $$ w_{1}(t)=t(\sigma[\rho]-[\rho u]),\,\,  \mathrm{and} \,\,\, w_{2}(t)=t(\sigma[\rho u]-[\rho u^{2}]),$$
 respectively, which form a delta shock solution of (1.5) with the same Riemann data (1.3).
Here $\sigma=\frac{\sqrt{\rho_{-}}u_{-}+\sqrt{\rho_{+}}u_{+}}{\sqrt{\rho_{-}}+\sqrt{\rho_{+}}}.$

\vskip 0.1in
\noindent{\small {\small\bf Proof.} (1) Set $\xi=\frac{x}{t}.$  Then for any fixed $\gamma \in(1,3)$, the Riemann solution
containing two shocks $S_{1}$ and $S_{2}$ of (1.7)  with the  Riemann initial data (1.3) can be written as $$(\rho_{\gamma},u_{\gamma})(\xi)=\left\{\begin{array}{ll} (\rho_{-},
u_{-}),\,\,\,\,\xi< \overline{\sigma}_{1},\\(\rho_{\ast},
u_{\ast}),\,\,\,\,\overline{\sigma}_{1}< \xi < \overline{\sigma}_{2},\\(\rho_{+},
u_{+}),\,\,\,\,\xi> \overline{\sigma}_{2}.\end{array} \right.\eqno{(6.18)}$$
From (5.2), we have the following weak formulations:
$$\int_{-\infty}^{+\infty}\rho_{\gamma}(\xi)(u_{\gamma}(\xi)-\xi)\varphi'(\xi)d\xi
-\int_{-\infty}^{+\infty}\rho_{\gamma}(\xi)\varphi(\xi)d\xi=0,
\eqno{(6.19)}$$
$$\int_{-\infty}^{+\infty}\rho_{\gamma}(\xi)u_{\gamma}(\xi)(u_{\gamma}(\xi)-\xi)\varphi'(\xi)d\xi
+\int_{-\infty}^{+\infty}(\rho_{\gamma}(\xi))^{\gamma}(u_{\gamma}(\xi)-\frac{1}{\gamma}\xi)\varphi'(\xi)d\xi
$$$$
-\int_{-\infty}^{+\infty}\left(\rho_{\gamma}(\xi)u_{\gamma}(\xi)+\frac{1}{\gamma}(\rho_{\gamma}(\xi))^{\gamma}\right)\varphi(\xi)d\xi=0,\eqno{(6.20)}
$$
for any $\varphi(\xi)\in C_{0}^{+\infty}(R)$.

(2) For the first integral on the left-hand side of (6.19),  using the method of integration by parts, we can derive
$$\int_{-\infty}^{+\infty}\rho_{\gamma}(\xi)(u_{\gamma}(\xi)-\xi)\varphi'(\xi)d\xi
=\left(\int_{-\infty}^{\overline{\sigma}_{1}}+\int_{\overline{\sigma}_{2}}^{+\infty}+\int_{\overline{\sigma}_{1}}^{\overline{\sigma}_{2}}\right)
\rho_{\gamma}(\xi)(u_{\gamma}(\xi)-\xi)\varphi'(\xi)d\xi
$$
$$=\rho_{-}u_{-}\varphi(\overline{\sigma}_{1})-\rho_{+}u_{+}\varphi(\overline{\sigma}_{2})
-\rho_{-}\overline{\sigma}_{1}\varphi(\overline{\sigma}_{1})+\rho_{+}\overline{\sigma}_{2}\varphi(\overline{\sigma}_{2})+\int_{-\infty}^{\overline{\sigma}_{1}}\rho_{-}\varphi(\xi)d\xi
$$
$$+\int_{\overline{\sigma}_{2}}^{+\infty}\rho_{+}\varphi(\xi)d\xi
+\int_{\overline{\sigma}_{1}}^{\overline{\sigma}_{2}}\rho_{\ast}(u_{\ast}-\xi)\varphi'(\xi)d\xi
$$
Meanwhile, we have
$$ \int_{\overline{\sigma}_{1}}^{\overline{\sigma}_{2}}\rho_{\ast}(u_{\ast}-\xi)\varphi'(\xi)d\xi
=\rho_{\ast}u_{\ast}(\varphi(\overline{\sigma}_{2})-\varphi(\overline{\sigma}_{1}))-\rho_{\ast}(\overline{\sigma}_{2}\varphi(\overline{\sigma}_{2})-\overline{\sigma}_{1}\varphi(\overline{\sigma}_{1}))
+\int_{\overline{\sigma}_{1}}^{\overline{\sigma}_{2}}\rho_{\ast}\varphi(\xi)d\xi
$$
$$=\rho_{\ast}(\overline{\sigma}_{2}-\overline{\sigma}_{1})\left(u_{\ast}\frac{\varphi(\overline{\sigma}_{2})-\varphi(\overline{\sigma}_{1})}{\overline{\sigma}_{2}-\overline{\sigma}_{1}}
+\frac{\int_{\overline{\sigma}_{1}}^{\overline{\sigma}_{2}}\varphi(\xi)d\xi}{\overline{\sigma}_{2}-\overline{\sigma}_{1}}
-\frac{\overline{\sigma}_{2}\varphi(\overline{\sigma}_{2})-\overline{\sigma}_{1}\varphi(\overline{\sigma}_{1})}{\overline{\sigma}_{2}-\overline{\sigma}_{1}}\right).
$$
Then, by Lemma 6.2-6,3, we can obtain
$$\lim_{\gamma\rightarrow1}\int_{\overline{\sigma}_{1}}^{\overline{\sigma}_{2}}\rho_{\ast}(u_{\ast}-\xi)\varphi'(\xi)d\xi=0.
$$
 Hence taking the limit $ \gamma\rightarrow1$ in (6.19) leads to
 $$\lim_{\gamma\rightarrow1}\int_{-\infty}^{+\infty}(\rho_{\gamma}(\xi)-\rho_{0}(\xi))\varphi(\xi)d\xi
=(\sigma[\rho]-[\rho u])\varphi(\sigma),\eqno (6.21)$$
where $ (\rho_{0}(\xi),u_{0}(\xi))=(\rho_{\pm},u_{\pm}),~\pm(\xi-\sigma)>0.$

(3) Similarly, we can obtain for (6.20) that
$$\lim_{\gamma\rightarrow1}
\int_{-\infty}^{+\infty}\rho_{\gamma}(\xi)u_{\gamma}(\xi)(u_{\gamma}(\xi)-\xi)\varphi'(\xi)d\xi
$$
$$=\left(\sigma[\rho u]-[\rho u^{2}]\right)\varphi(\sigma)+\int_{-\infty}^{+\infty}\rho_{0}(\xi)u_{0}(\xi)\varphi(\xi)d\xi,
$$
and
$$\int_{-\infty}^{+\infty}(\rho_{\gamma}(\xi))^{\gamma}(u_{\gamma}(\xi)-\frac{1}{\gamma}\xi)\varphi'(\xi)d\xi
=\left(\int_{-\infty}^{\overline{\sigma}_{1}}+\int_{\overline{\sigma}_{2}}^{+\infty}+\int_{\overline{\sigma}_{1}}^{\overline{\sigma}_{2}}\right)
(\rho_{\gamma}(\xi))^{\gamma}\left(u_{\gamma}(\xi)-\frac{1}{\gamma}\xi\right)\varphi'(\xi)d\xi
$$
$$=\rho^{\gamma}_{-}u_{-}\varphi(\overline{\sigma}_{1})-\rho^{\gamma}_{+}u_{+}\varphi(\overline{\sigma}_{2})
-\frac{1}{\gamma}\rho^{\gamma}_{-}\overline{\sigma}_{1}\varphi(\overline{\sigma}_{1})+\frac{1}{\gamma}\rho^{\gamma}_{+}\overline{\sigma}_{2}\varphi(\overline{\sigma}_{2})
+\int_{-\infty}^{\overline{\sigma}_{1}}\frac{1}{\gamma}\rho^{\gamma}_{-}\varphi(\xi)d\xi
$$
$$ +\int_{\overline{\sigma}_{2}}^{+\infty}\frac{1}{\gamma}\rho^{\gamma}_{+}\varphi(\xi)d\xi
+\frac{1}{\gamma}\rho^{\gamma}_{\ast}(\overline{\sigma}_{2}-\overline{\sigma}_{1})
\left(\gamma u_{\ast}\frac{\varphi(\overline{\sigma}_{2})-\varphi(\overline{\sigma}_{1})}{\overline{\sigma}_{2}-\overline{\sigma}_{1}}
-\frac{\overline{\sigma}_{2}\varphi(\overline{\sigma}_{2})-\overline{\sigma}_{1}\varphi(\overline{\sigma}_{1})}{\overline{\sigma}_{2}-\overline{\sigma}_{1}}
+\frac{\int^{\overline{\sigma}_{2}}_{\overline{\sigma}_{1}}\varphi(\xi)d\xi}{\overline{\sigma}_{2}-\overline{\sigma}_{1}}\right),
$$
which converges to
$$(\sigma[\rho]-[\rho u])\varphi(\sigma)+\int_{-\infty}^{+\infty}\rho_{0}(\xi)\varphi(\xi)d\xi
$$
by Lemma 6.2-6.3.

Thus, from (6.21), we can get
$$ \lim\limits_{\gamma\rightarrow1}\int_{-\infty}^{+\infty}(\rho_{\gamma}(\xi)u_{\gamma}(\xi)-\rho_{0}(\xi)u_{0}(\xi))\varphi(\xi)d\xi
=\left(\sigma[\rho u]-[\rho u^{2}]\right)\varphi(\sigma).\eqno (6.22)
$$

(4) Finally, we study the limits of $\rho_{\gamma}(t,x) $ and $\rho_{\gamma}(t,x)u_{\gamma}(t,x) $ depending on $t$ as $\gamma\rightarrow1$.
Regarding $t$ as a parameter, we can  get from   (6.21) that
$$
\lim_{\gamma\rightarrow1}\int_{-\infty}^{+\infty}(\rho_{\gamma}(\xi)-\rho_{0}(\xi))\varphi(t,\xi t)d\xi
=\lim_{\gamma\rightarrow1}\int_{-\infty}^{+\infty}(\rho_{\gamma}(x/t)-\rho_{0}(x/t))\varphi(t,x)d(x/t)
$$
$$=\frac{1}{t}\lim_{\gamma\rightarrow1}\int_{-\infty}^{+\infty}(\rho_{\gamma}(t,x)-\rho_{0}(t,x))\varphi(t,x)dx
=(\sigma[\rho]-[\rho u])\varphi(t, \sigma t).\eqno (6.23)
$$
Then multiplying (6.23) by $t$ and taking integration, we have
$$\lim_{\gamma\rightarrow1}\int_{0}^{+\infty}\int_{-\infty}^{+\infty}(\rho_{\gamma}(t,x)-\rho_{0}(t,x))\varphi(t,x)dx dt
=\int_{0}^{+\infty}t(\sigma[\rho]-[\rho u])\varphi(t,\sigma t)dt
$$
in which by definition (2.3), we have
$$\int_{0}^{+\infty}t(\sigma[\rho]-[\rho u])\varphi(t,\sigma t)dt=
\langle w_{1}(\cdot)\delta_{S},\varphi(\cdot,\cdot)\rangle.\eqno (6.24)
$$
where
 $$ w_{1}(t)=t(\sigma[\rho]-[\rho u]).$$

In the same way, we can derive from (6.22) that
$$\lim_{\gamma\rightarrow1}\int_{0}^{+\infty}\int_{-\infty}^{+\infty}(\rho_{\gamma}(t,x)u_{\gamma}(t,x)-\rho_{0}u_{0}(t,x))\varphi(t,x)dxdt
=\langle w_{2}(\cdot)\delta_{S},\varphi(\cdot,\cdot)\rangle.\eqno (6.25)
$$
where
 $$ w_{2}(t)=t(\sigma[\rho u]-[\rho u^{2}]).$$
The proof is completed. $~~\Box$

\baselineskip 15pt
 \sec{\Large\bf 7.\quad  Numerical results}
  In this section, we use the fifth-order weighted essentially non-oscillatory scheme and third-order Runge-Kutta method  [12, 27]  with the mesh 400 points to
present   some  groups of representative numerical results
  for the Aw-Rascle traffic model (1.1)-(1.2) and the perturbed Aw-Rascle model (1.7) as $\gamma$ decreases. A number of iterative numerical trials are executed to guarantee what we demonstrate are not numerical objects. The numerical simulations are consistent with the theoretical analysis.

\baselineskip 15pt
 \sec{\Large\bf 7.1.\quad Formation of delta-shocks  in (1.1)-(1.2)}

The numerical simulations are corresponding to the theoretical analysis in Section 4.
When $(\rho_{+},u_{+})\in I(\rho_{-},u_{-})$, we take the initial data as follows:
$$ (\rho, u)(0, x) =\left\{\begin{array}{ll} (3.5,
6),\,\,\,\,\,x< 0,\\(2,
4),\,\,\,\,x> 0,\end{array} \right.\eqno{(7.1)}$$
 and compute the solution of the  Riemann problem of (1.1)-(1.2) up to $t=0.4$, the numerical simulations for different choices of  $\gamma$, starting with $\gamma$=0.6,  then  $\gamma$=0.3, and finally $\gamma$=0.01,  are presented in Figs. 3-5 which show the process of concentration  and formation of the delta shock wave  in vanishing adiabatic exponent limit of solutions   containing a shock wave and a contact discontinuity.

\begin{figure}[htbp]
\centering
\begin{minipage}[c]{0.33\textwidth}
\centering
\includegraphics[width=2in]{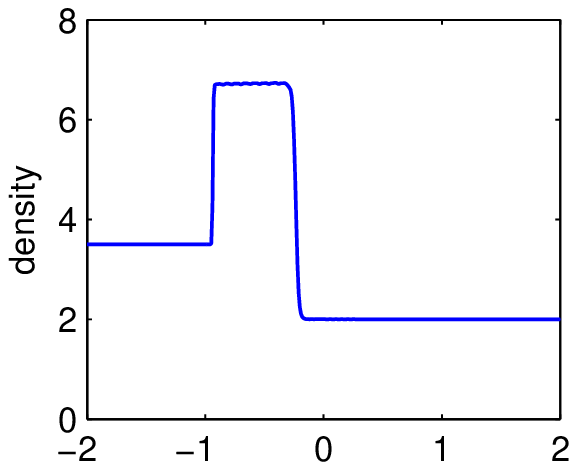}
\end{minipage}%
\begin{minipage}[c]{0.33\textwidth}
\centering
\includegraphics[width=2in]{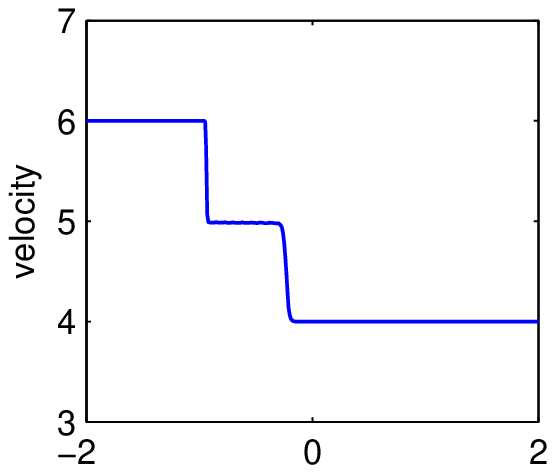}
\end{minipage}%
\begin{minipage}[c]{0.28\textwidth}
\centering
\includegraphics[width=2in]{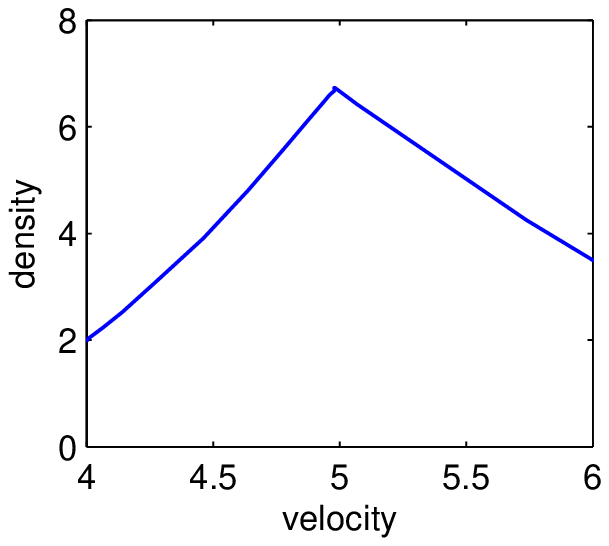}
\end{minipage}%
\end{figure}

\centerline{\bf Fig. 3.\,\,   Density (left) and velocity (right) for $\gamma=0.6$.}

\begin{figure}[htbp]
\centering
\begin{minipage}[c]{0.33\textwidth}
\centering
\includegraphics[width=2in]{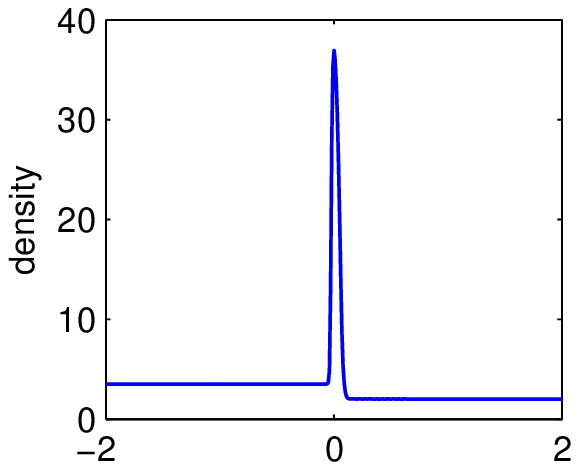}
\end{minipage}%
\begin{minipage}[c]{0.33\textwidth}
\centering
\includegraphics[width=2in]{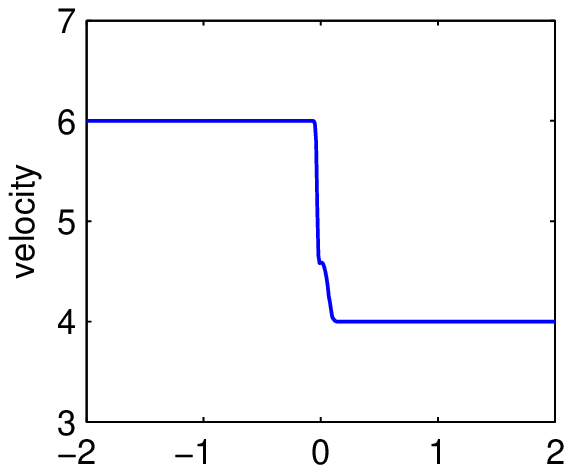}
\end{minipage}%
\begin{minipage}[c]{0.28\textwidth}
\centering
\includegraphics[width=2in]{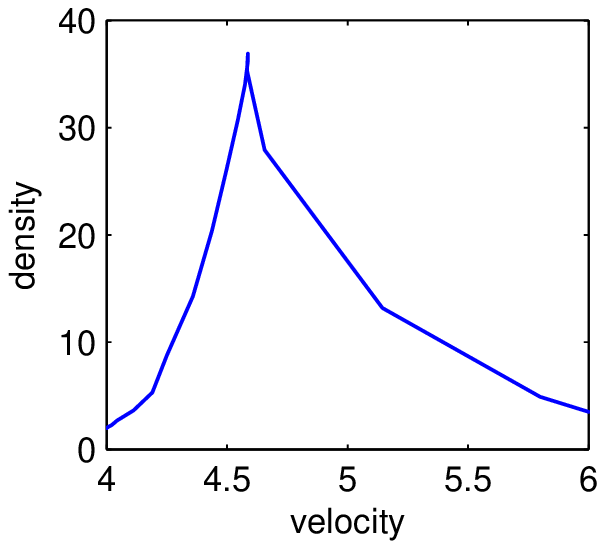}
\end{minipage}%
\end{figure}

\centerline{\bf Fig. 4.\,\,   Density (left) and velocity (right) for $\gamma=0.3$.}

\newpage

\begin{figure}[htbp]
\centering
\begin{minipage}[c]{0.33\textwidth}
\centering
\includegraphics[width=2in]{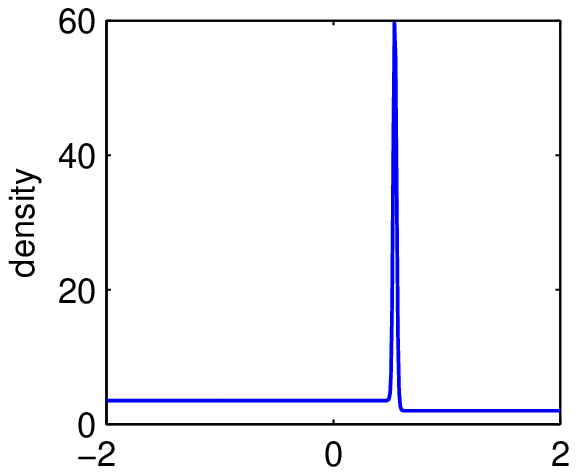}
\end{minipage}%
\begin{minipage}[c]{0.33\textwidth}
\centering
\includegraphics[width=2in]{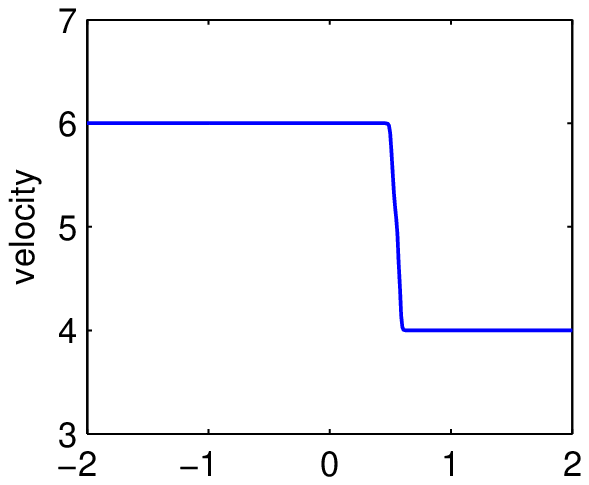}
\end{minipage}%
\begin{minipage}[c]{0.28\textwidth}
\centering
\includegraphics[width=2in]{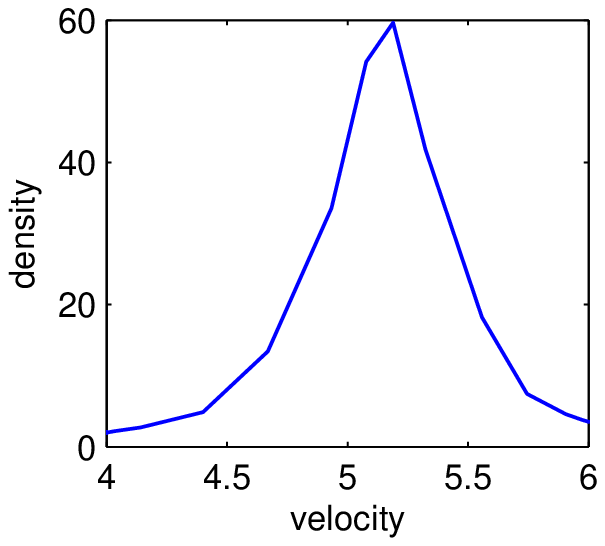}
\end{minipage}%
\end{figure}

\centerline{\bf Fig. 5.\,\,   Density (left) and velocity (right) for $\gamma=0.001$.}

From these numerical results, we can clearly observe that,  when $\gamma$ decreases, the locations of the
shock wave and contact discontinuity become closer and closer, and the density of the intermediate state
increases dramatically, while  the velocity becomes a piecewise constant function. In the end, as $\gamma\rightarrow0$, along with the intermediate state, the shock wave and the contact discontinuity coincide to form a delta-shock, while the velocity keeps a step function. The numerical simulations are in complete agreement with the theoretical analysis in Section 4.

\baselineskip 15pt
 \sec{\Large\bf 7.2.\quad Formation of delta-shocks  in (1.7)}
The numerical simulations are corresponding to the theoretical analysis in Section 6.
When $(\rho_{+},u_{+})\in S_{1}S_{2}(\rho_{-},u_{-})$, we take the initial data as follows:
$$ (\rho, u)(0, x) =\left\{\begin{array}{ll} (3,
4),\,\,\,\,\,x< 0,\\(2.5,
2),\,\,\,x> 0,\end{array} \right.\eqno{(7.2)}$$
 and compute the solution of the  Riemann problem of (1.7) up to $t=0.4$, the numerical simulations for different choices of  $\gamma$, starting with $\gamma$=1.4, then $\gamma$= 1.04, and finally  $\gamma$= 1.001, are presented in Figs. 6-8 which show the process of concentration  and formation of the delta shock wave  in the pressureless limit of solutions containing two shocks.

\begin{figure}[htbp]
\centering
\begin{minipage}[c]{0.33\textwidth}
\centering
\includegraphics[width=2in]{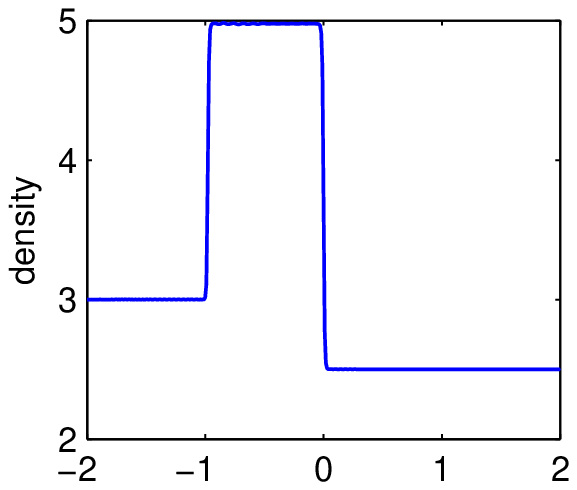}
\end{minipage}%
\begin{minipage}[c]{0.33\textwidth}
\centering
\includegraphics[width=2in]{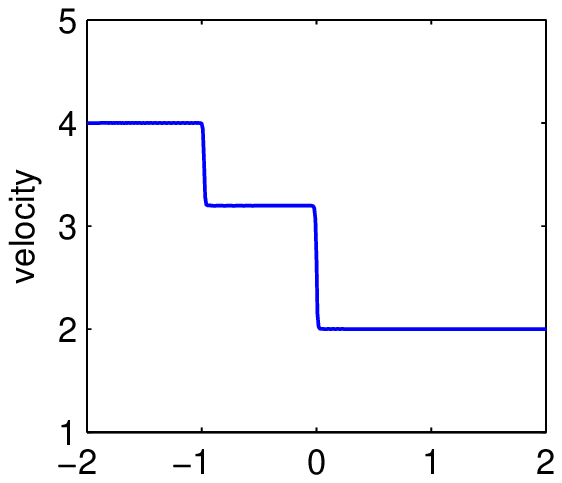}
\end{minipage}%
\begin{minipage}[c]{0.28\textwidth}
\centering
\includegraphics[width=2in]{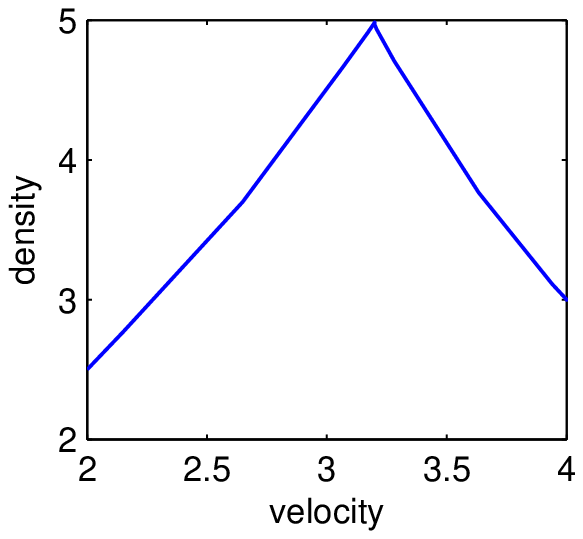}
\end{minipage}%
\end{figure}

\centerline{\bf Fig. 6.\,\,   Density (left) and velocity (right) for $\gamma=1.4$.}

\newpage
\begin{figure}[htbp]
\centering
\begin{minipage}[c]{0.33\textwidth}
\centering
\includegraphics[width=2in]{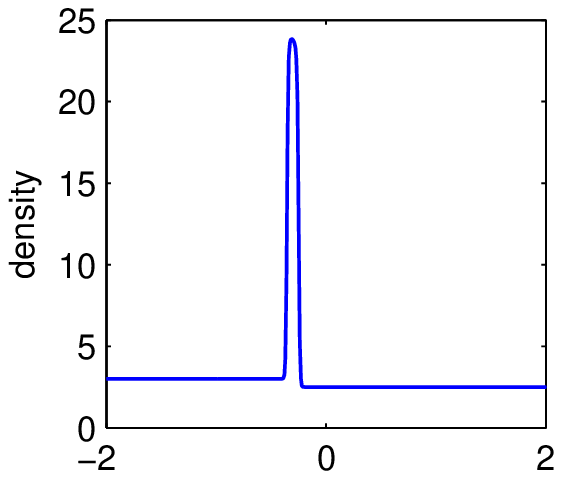}
\end{minipage}%
\begin{minipage}[c]{0.33\textwidth}
\centering
\includegraphics[width=2in]{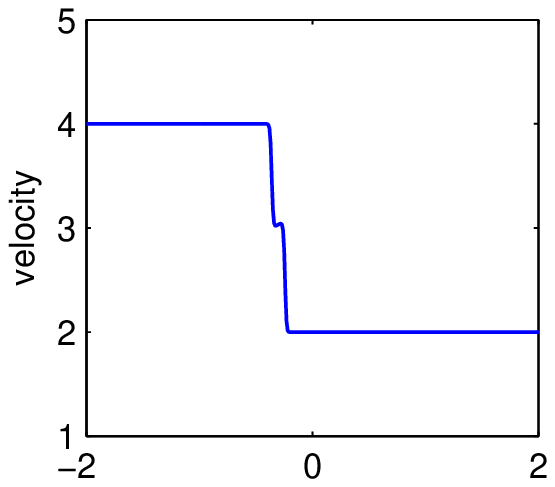}
\end{minipage}%
\begin{minipage}[c]{0.28\textwidth}
\centering
\includegraphics[width=2in]{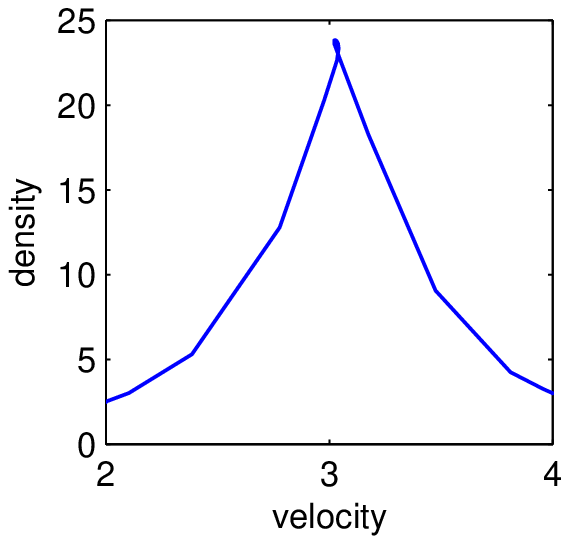}
\end{minipage}%
\end{figure}

\centerline{\bf Fig. 7.\,\,   Density (left) and velocity (right) for $\gamma=1.04$.}

\begin{figure}[htbp]
\centering
\begin{minipage}[c]{0.33\textwidth}
\centering
\includegraphics[width=2in]{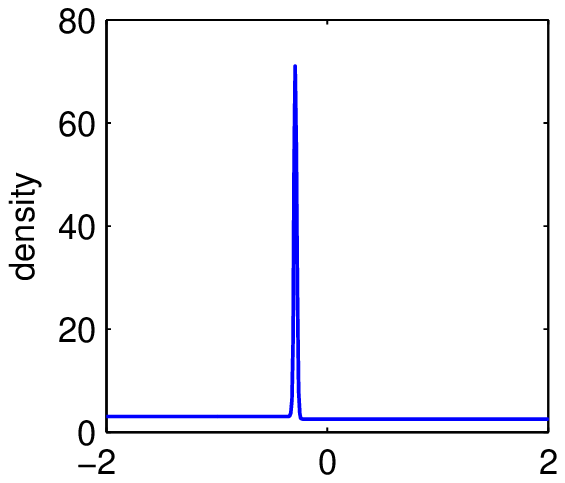}
\end{minipage}%
\begin{minipage}[c]{0.33\textwidth}
\centering
\includegraphics[width=2.0in]{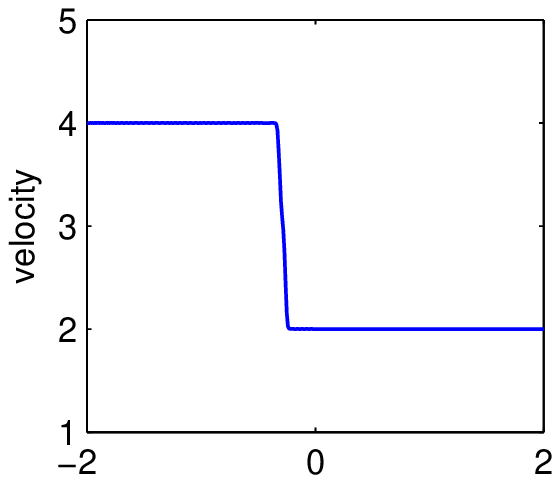}
\end{minipage}%
\begin{minipage}[c]{0.28\textwidth}
\centering
\includegraphics[width=2.1in]{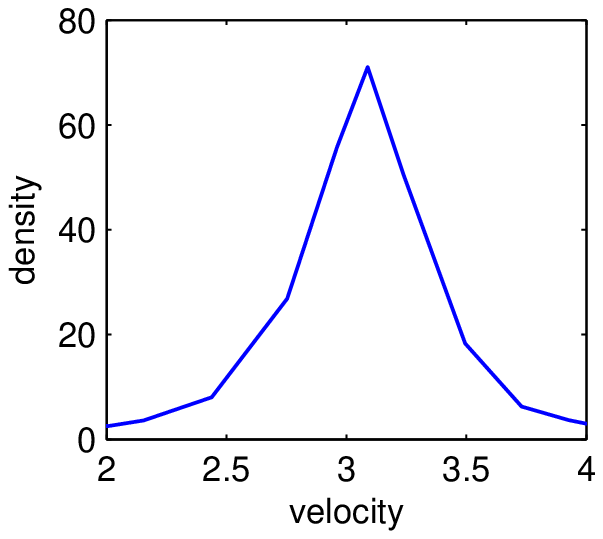}
\end{minipage}%
\end{figure}

\centerline{\bf Fig. 8.\,\,   Density (left) and velocity (right) for $\gamma=1.001$.}

From these numerical results, we can clearly observe that, as $\gamma$ decreases, the locations of the two shocks become closer and closer, and the density of the intermediate state increases dramatically, while  the velocity becomes a piecewise constant function. In the end, as $\gamma\rightarrow1$, along with the intermediate state, the two shocks coincide to form  the delta shock wave of  the zero pressure gas dynamics 
 (1.5), while the velocity keeps a step function. The numerical simulations are in complete agreement with the theoretical analysis in Section 6.


\newpage
  \vskip 10 pt
\begin{thebibliography}{99}

\bibitem{s3}

 A. Aw, M. Rascle, Resurrection of ``second order" models of traffic flow, SIAM J. Appl. Math. 60 (2000) 916-938.
\bibitem{s1}  F. Bouchut, On zero pressure gas dynamics, in: Advances in Kinetic Theory and Computing, in: Ser. Adv. Math. Appl. Sci.,
vol. 22, World Scientific Publishing, River Edge, NJ, 1994, pp. 171-190.


\bibitem{s1} Y. Brenier, E. Grenier, Sticky particles and scalar conservation laws, SIAM J. Numer. Anal. 35 (1998) 2317-2328.

\bibitem{s7}
 G.Q. Chen,  H. Liu, Formation of ${\delta}$-shocks and vacuum states in the vanishing pressure limit of solutions to the Euler equations for isentropic fluids,  SIAM J. Math. Anal.  34 (2003) 925-938.

 \bibitem{A6}
G.-Q. Chen,  H. Liu,
 Concentration and cavitation in the vanishing pressure limit of solutions to the Euler equations for nonisentropic fluids, Phys.
D 189  (2004) 141-165.
\bibitem{s1} C. Daganzo, Requiem for second order fluid approximations of traffic flow, Transportation Res. Part B 29 (1995) 277-286.
\bibitem{s3}W. E, Yu.G. Rykov, Ya.G. Sinai, Generalized variational principles, global weak solutions and behavior with random initial data for systems of conservation laws arising in adhesion particle dynamics, Comm. Math. Phys. 177 (1996) 349-380.

\bibitem{s3} S. Ha, F. Huang, and Y. Wang, A global unique solvability of entropic weak solution to the one-dimensional pressureless
Euler system with a flocking dissipation, J. Differ. Equations 257 (2014) 1333-1371 .

\bibitem{s3} F. Huang,  Z. Wang, Well posedness for pressureless flow, Comm. Math. Phys. 222 (2001)
 117-146.

\bibitem{s1} M. Ibrahim, F. Liu, S. Liu,
Concentration of mass in the pressureless limit of Euler equations for power law,
Nonlinear Anal. Real World Appl.  47 (2019) 224-235.


\bibitem{s3} K.T. Joseph, A Riemann problem whose viscosity solutions contain $\delta$-measures, Asymptot Anal.  7(1993) 105-120.
\bibitem{s1}A. Kurganov, E. Tadmor, New high-resolution central schemes for nonlinear conservation laws and
convection diffusion equations, J. Comput. Phys. 160 (2000) 241-282.



\bibitem{s1} J. Lebacque, S. Mammar, and H. Salem, The Aw-Rascle and Zhang¡¯s model: Vacuum problems, existence and regularity
of the solutions of the Riemann problem, Transp. Res. Part B 41 (2007) 710-721.


\bibitem{s1}J. Li, Note on the compressible Euler equations with zero temperature, Appl. Math. Lett.
14 (2001) 519-523.
\bibitem{s1} J. Li, H. Yang, Delta-shocks as limits of vanishing viscosity for multidimensional zero-pressure
gas dynamics, Quart. Appl. Math. 59 (2) (2001) 315-342.

\bibitem{s1} J. Li, T. Zhang, S. Yang, The two-dimensional Riemann problem in gas dynamics, Vol. 98
of Pitman Monographs and Surveys in Pure and Applied Mathematics, Longman, Harlow,
1998.

\bibitem{s1} H. Li, Z. Shao, Delta shocks and vacuum states in vanishing pressure limits of solutions to the
relativistic Euler equations for generalized Chaplygin gas, Commun. Pure Appl. Anal. 15 (2016)
2373-2400.
\bibitem{s1}
J. Liu,  W. Xiao,
Flux approximation to the Aw-Rascle model of traffic flow,
Journal of Mathematical Physics 59, 101508 (2018); doi: 10.1063/1.5063469.
\bibitem{s1} D. Mitrovic, M. Nedeljkov, Delta-shock waves as a limit of shock waves, J. Hyperbolic Differ. Equ. 4 (2007) 629-653.
\bibitem{s3}
L. Pan, X. Han,
The Aw-Rascle traffic model with Chaplygin pressure,
J. Math. Anal. Appl. 401 (2013) 379-387.




\bibitem{s1}S.F. Shandarin, Ya.B. Zeldovich, The large-scale structure of the universe: turbulence, intermittency, structures in a self-gravitating medium, Rev. Modern Phys. 61 (1989) 185-220.
















\bibitem{s1} C. Shen, The limits of Riemann solutions to the isentropic magnetogasdynamics, Appl. Math. Lett. 24 (2011) 1124-1129.


\bibitem{s1} C. Shen, M. Sun, Formation of delta shocks and vacuum states in the vanishing pressure limit of Riemann solutions to the perturbed Aw-Rascle model, J. Differential Equations 249 (2010) 3024-3051.
  \bibitem{s3} C. Shen, M. Sun,  Z. Wang, Limit relations for three simple hyperbolic systems of
conservation laws, Math. Meth. Appl. Sci. 33 (2010) 1317-1330.



\bibitem{s1} W. Sheng, G. Wang, G. Yin, Delta wave and vacuum state for generalized Chaplygin gas dynamics system as pressure vanishes, Nonlinear Anal. Real World Appl. 22 (2015) 115-128.



\bibitem{s3}
 W. Sheng, T. Zhang, The Riemann problem for the transportation equations in gas dynamics, in: Mem. Amer. Math. Soc., 137, AMS,
Providence, 1999.
\bibitem{s1}C. W. Shu, Essentially non-oscillatory and weighted essentially non-oscillatory schemes for hyperbolic
conservation laws, in Advanced Numerical Approximation of Nonlinear  Hyperbolic Equations, Lecture Notes in Mathematics
 Vol. 1697 (Springer Berlin Heidelberg, 1998), pp. 325-432.




\bibitem{s1}
M. Sun,  Interactions of elementary waves for the Aw-Rascle model,
SIAM J. Appl. Math.
 69 (2009) 1542-1558.


\bibitem{s1} G. Yin, W. Sheng, Delta shocks and vacuum states in vanishing pressure limits of solutions to the relativistic Euler
equations for polytropic gases, J. Math. Anal. Appl. 355  (2009) 594-605.

\bibitem{s1} H. Zhang, A non-equilibrium traffic model devoid of gas-like behavior, Transportation Res. Part B 36 (2002) 275-290.
















\end{thebibliography}
\end{document}